\documentclass[12pt]{iopart}

\usepackage{iopams}  
\usepackage{subfigure}
\usepackage[T1]{fontenc}
\usepackage{dsfont,graphicx}
\newtheorem{proposition}{Proposition}[section]
\newtheorem{theorem}{Theorem}[section]
\newtheorem{lemma}{Lemma}

\newtheorem{dfn}{Definition}

\newcommand{\deq}{\mathrel{\mathop:} = } 
\newcommand{\paren}[1]{{\left( #1 \right)}} %

\newcommand{\N}{\mathbb{N}}

\begin{document}

\title[Adaptive estimation of density matrix in QHT]{Adaptive estimation of the
density matrix in quantum homodyne tomography with noisy data}

\author{P Alquier$^1$, K Meziani$^2$ and G Peyr\'e$^3$}

\address{$^1$ School of Mathematical Sciences, University College, Dublin.}
\address{ $^2$$^3$ CEREMADE, UMR CNRS 7534, Universit\'e Paris Dauphine, Place du Mar\'echal De Lattre De Tassigny, 75775 PARIS Cedex 16, France.
}
\eads{\mailto{pierre.alquier@ucd.ie}, \mailto{meziani@ceremade.dauphine.fr}, \\
\mailto{gabriel.peyre@ceremade.dauphine.fr}}
\begin{abstract}
     In the framework of noisy quantum homodyne tomography with
    efficiency parameter $1/2 < \eta \leq 1$, we propose a novel estimator
    of a quantum state whose density matrix elements $\rho_{m,n}$
    decrease like $Ce^{-B(m+n)^{ r/ 2}}$, for fixed $C\geq 1$, $B>0$ and
    $0<r\leq 2$. On the contrary to previous works, we focus on the
    case where $r$, $C$ and $B$ are unknown. The procedure estimates
    the matrix coefficients by a projection method on the pattern functions,
    and then by soft-thresholding the estimated coefficients. 
    We prove that under the  $\mathbb{L}_2$ -loss our procedure is adaptive rate-optimal, in the sense that it achieves the same rate of conversgence as the best possible procedure relying on the knowledge of $(r,B,C)$. Finite sample behaviour of our adaptive procedure are explored through numerical experiments.
\end{abstract}


\vspace{2pc}


\maketitle

\section{Introduction}

This paper deals with a \textit{severely ill-posed inverse problem} which comes from
quantum optics. 
Quantum optics is a branch of quantum mechanics which studies  physical systems
at the atomic and subatomic scales. As the language used by physicists\footnote{(e.g. they speak about `states'' or ``observable'' instead of
``laws'' or  ``random variables''...)} differs
from the one that is used by statisticians, we start with general notions on
quantum mechanics. 
The interested reader can get further acquaintance with quantum concepts
through the textbooks  or the review articles
\cite{Helstrom,Holevo,Barndorff-Nielsen&Gill&Jupp,Leonhardt}.

\subsection{Physical background}
\label{physbackground}

\subsubsection{Quantum mechanics}
\label{Quantum mechanics}

In quantum mechanics, the quantum state of a system is a mathematical object 
which encompasses all the information about the system.  The most common
representation of a quantum state is an operator $\rho$ on a complex Hilbert
space $\mathcal{H}$ (called the space of states) satisfying the three following conditions:
\begin{enumerate}
     \item Self adjoint: $\rho=\rho^{*}$, where $\rho^{*}$ is the adjoint of
$\rho$.
     \item Positive: $\rho\geq 0$, or equivalently
$\langle\psi,\rho\psi\rangle\geq0$ for all 
               $\psi\in \mathcal  {H}$.
     \item Trace one: $\mathrm{Tr}(\rho)=1$.
\end{enumerate}
A quantum state $\rho$ encodes the probabilities of the measurable properties,
or ``\textbf{observables}'' (energy, position, ...) of the considered quantum
system. Generally, in quantum mechanics the expected results of the measurements
of an observable are not  deterministic values but  predictions about
probability distributions, that is the probability of obtaining each of the
possible outcomes when measuring an observable. \\
An observable $\mathbf{X}$ is described by a self adjoint operator on the space
of states $\mathcal{H}$ and
$$
\mathbf{X}=\sum_{a}^{dim\mathcal{H}}x_a\mathbf{P}_a,
$$
where the eigenvalues $\{x_a\}_a$ of the observable $\mathbf{X}$ are real and 
$\mathbf{P}_a$  is the projection onto the one dimensional space generated by
the eigenvector of $\mathbf{X}$ corresponding to the eigenvalue $x_a$. Then, when performing a measurement of the
observable $\mathbf{X}$ of a quantum state $\rho$, the result is a random
variable $X$ with values in the set of the eigenvalues of the observable
$\mathbf{X}$. For a quantum system prepared in state $\rho$,  $X$ has the
following probability distribution and expectation  function
\begin{center}
$\mathbb{P}_\rho(X=x_a)=\mathrm{Tr}(\mathbf{P}_a\rho)\quad$ and
$\quad\mathbb{E}_\rho(X)=\mathrm{Tr}(\mathbf{X}\rho)$.
\end{center}
\noindent\\
 An important element which affects the result of the measurement process is the
purity of quantum states. A state is called pure if it cannot be represented as
a mixture (convex combination) of other states, i.e., if it is an extreme point
of the convex set of states. All other states are called mixed states. We give
examples of states in Section~\ref{sec.simulations}.

\subsubsection{Quantum optics}
\label{Quantum opticcs}

In this paper, the quantum system we work with is a monochromatic light in a
cavity described by a quantum harmonic oscillator. In the framework of quantum
optics,  the space of states is known to be the separable Hilbert space
$\mathcal{H}=\mathbb{L}_2(\mathbb{R})$, \textit{i.e.} the space of square integrable
complex valued functions on the real line. A particular orthonormal basis
$\left(\psi_j\right)_{j\in\mathbb{N}}$ -- called the Fock basis -- comes with this Hilbert space. This
physically very meaningful basis is defined for all $j\in\mathbb{N}$ as follows
          \begin{equation}
          \label{eq:fock}
          \psi_j(x)\deq \frac{1}{\sqrt{\sqrt{\pi}2^jj!}}H_j(x)e^{-x^2/2},
          \end{equation}
where $H_j(x) \deq (-1)^j e^{x^2} \frac{d^j}{dx^j} e^{-x^2}$ is the $j$-th
Hermite polynomial. In the Fock basis~(\ref{eq:fock}), a state is described  by an
infinite density matrix $\rho=[\rho_{j,k}]_{j,k\in\N}$. \\

\noindent We may give  an equivalent  representation for a quantum state $\rho$ in terms of the
associated Wigner function $W_\rho$ (see \cite{Wigner}). The Wigner function
$W_\rho$  is a real function of two variables and may be defined by its Fourier
transform $\mathcal{F}_2$ with respect to both variables
        \begin{eqnarray*}
         \label{def.Wigner}        
\widetilde{W}_{\rho}(u,v):=\mathcal{F}_2[W_\rho](u,v)=\mathrm{Tr}
\left(\rho\exp(iu\textbf{Q}+iv  \textbf{P})\right),
         \end{eqnarray*}
where $\textbf{Q}$ and $\textbf{P}$ are respectively the electric and magnetic
fields.  These two observables, we are concerned by, do not commute. As
non-commuting observables,  they may not be simultaneously measurable. Therefore,
by performing measurements on $(\textbf{Q},\textbf{P})$, we cannot get a probability density of the result $(Q,P)$.  However, for $\phi\in[0,\pi]$ we can measure
the quadrature observables $\mathbf{X}_\phi :=\mathbf{Q}\cos \phi
+\mathbf{P}\sin \phi$, and then the above Wigner function plays the role of a
quasi-probability density. It does not satisfy all the properties of a
conventional probability density but satisfies boundedness properties
unavailable for classical densities. For instance, the Wigner function
can and normally does go negative for states which have no classical model. The
Wigner function is such that
        \begin{itemize}
        \item $W_\rho:\mathbb{R}^2\rightarrow\mathbb{R}$
        \item $\int\int W_\rho(q,p)dqdp=1$,
        \end{itemize}
Furthermore, its Radon transform is always a probability density
        \begin{equation}
        \label{eq.Radon.transform}
         p_\rho(x|\phi) \deq
\mathcal{R}[W_{\rho}](x,\phi)=\int_{-\infty}^\infty W_{\rho}(
        x\cos\phi - t\sin\phi, \, x \sin\phi + t\cos\phi )dt,
         \end{equation}
with respect to $\frac{1}{\pi}\lambda$, $\lambda$ being the Lebesgue measure on
$\mathbb{R}\times [0,\pi]$. \\

\noindent Now we can make explicit the links between the state $\rho$ and  the Radon transform $p_\rho(x|\phi)$ of the Wigner function $W_{\rho}$ associated to $\rho$.  In the Fock
basis (\ref{eq:fock}),  the entries $\rho_{j,k}$ of the infinite density matrix $\rho$ are given by
       \begin{eqnarray}
       \label{rhojk}
       \rho_{j,k}=\frac 1\pi\int \int_0^\pi
p_\rho(x|\phi)f_{j,k}(x)e^{-i(k-j)\phi}d\phi dx
       \end{eqnarray}
for all $j,k \in\N$.  The functions $f_{j,k} = f_{k,j}$, in the expression
(\ref{rhojk}), are bounded real functions commonly called \textit{pattern
functions} in quantum homodyne literature. A concrete expression for their
Fourier transform $\tilde{f}_{k,j}$ using Laguerre polynomials $L_{n}^{\alpha}(\cdot)$ is ( cf \cite{Richter1}):
for $j \geq k$,
        \begin{eqnarray}
         \label{Rpattern}
         \tilde{f}_{k,j}(t)  &=&\pi(-i)^{j-k}\sqrt{\frac{2^{k-j} k!}{j!}}|t|
          t^{j-k}e^{-\frac{t^2}{4}}L^{j-k}_k(\frac{t^2}{2}).
         \end{eqnarray}
We recall that the Laguerre polynomial of degree $n$ and
order $\alpha$ is defined by
$$L_{n}^{\alpha}(x) \deq (n!)^{-1} e^x x^{-\alpha} \frac{d^n}{dx^n}
(e^{-x} x^{n+\alpha}).$$

\subsubsection{Quantum Homodyne Tomography}
\label{QHT}

 In this paper, we address the problem of reconstructing the density matrix
$\rho$  of a monochromatic light in a cavity.  As the observables $\mathbf{Q}$
and $\mathbf{P}$ cannot  be measured simultaneously,  we measure the quadrature 
$\mathbf{X}_\phi :=\mathbf{Q}\cos \phi +\mathbf{P}\sin \phi$, where
$\phi\in[0,\pi]$. Each of these quadratures could be measured on a laser beam by
a technique put in practice for the first time in \cite{Smithey} and called
Quantum Homodyne Tomography (QHT). The theoretical foundation of quantum
homodyne tomography was outlined in \cite{Vogel&Risken}.\\

\begin{figure}[!h]
\begin{center}
\includegraphics[width=10cm,height=6cm]{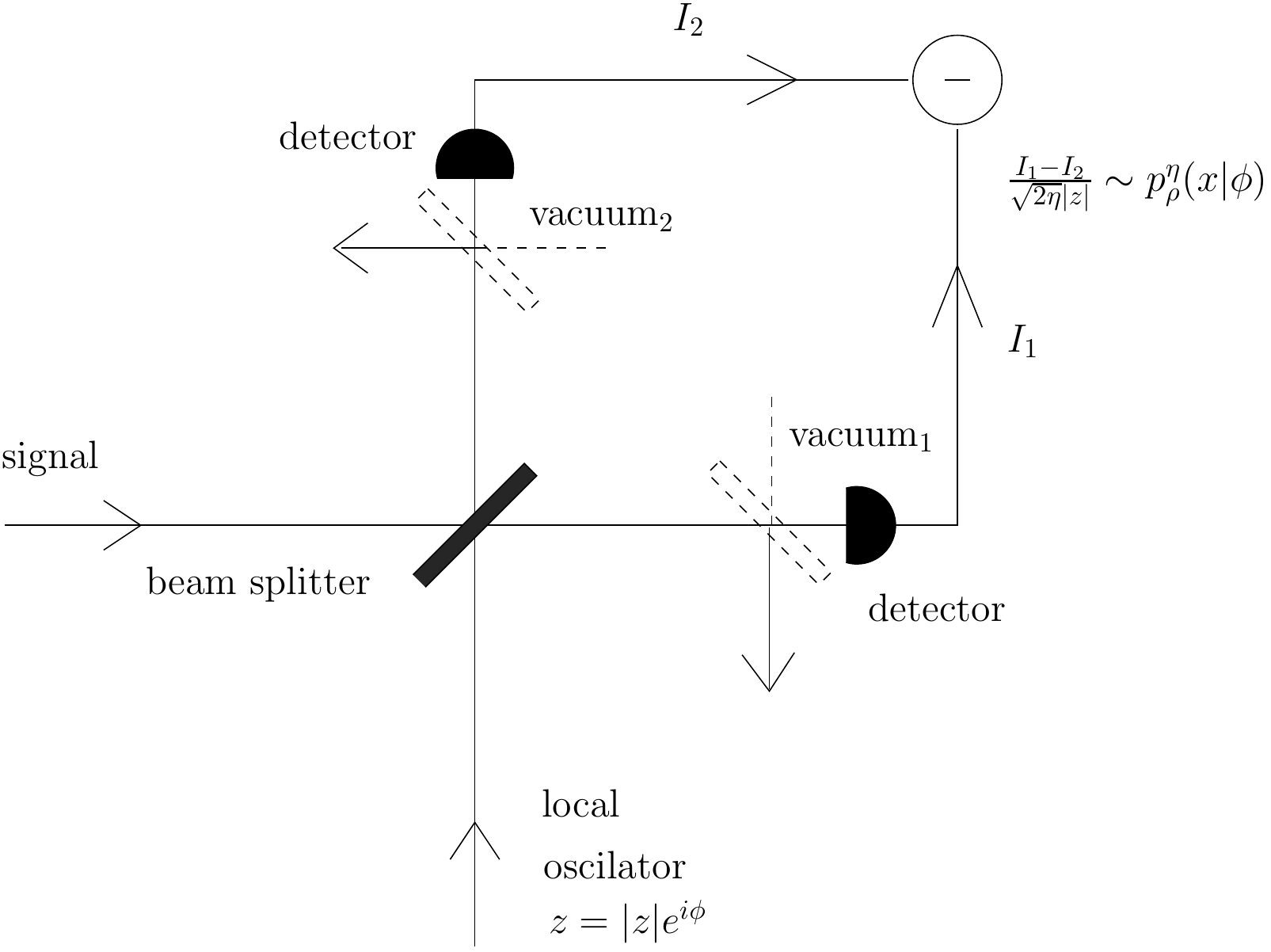}
\caption{QHT measurement scheme}
\label{Pfig:1}
\end{center}
\end{figure}

\noindent The experimental set-up, described in Figure~\ref{Pfig:1}, consists of mixing the
cavity pulse prepared in state $\rho$ with an additional laser of high intensity
$\left|z\right|>>1$ called the local oscillator.
 After the mixing, the beam is split again and each of the two emerging beams is
measured by  one of the two photodetectors which give integrated currents $I_1$ 
and $I_2$ proportional to the number of photons. The result of the measurement is
produced by taking the difference of the two currents and rescaling it by the
intensity $|z|$.  Just before the mixing the experimentalist may choose the
phase $\Phi$ of the local oscillator, randomly, uniformly distributed on
$[0,\pi]$.  In the case of noiseless measurement and for a phase $\Phi=\phi$,
the result $X_\phi=\frac{I_2-I_1}{|z|}$ has density $p_\rho(x|\phi)$
corresponding to measuring $\mathbf{X}_\phi$. \\

\noindent  In practice, a number of photons fails to be detected. These losses may be
quantified by one single coefficient $\eta\in[0,1]$, such that  $\eta=0$ when
there is no detection and  $\eta=1$  corresponds to the ideal case (no loss). The
physicists argue, that their machines actually have high detection efficiency,
around 0.8/0.9. Thus, we suppose $\eta$ known. As the detection process is
inefficient, an independent gaussian noise interferes additively with the ideal
data $X_\phi$. Thus for $\Phi=\phi$, the effective result of the QHT measurement
(Figure~\ref{Pfig:1}) is   for a known efficiency $\eta\in]0.5,1]$,  
$$
Y=\sqrt{\eta}\,X_\phi+\sqrt{(1-\eta)/2}\,\xi
$$
where $\xi$ is a standard Gaussian random variable, independent of $X_\phi$.

\subsection{Statistical model}
\label{Statistical.model}

This paper aims at  reconstructing the density matrix of  a monochromatic light
in a cavity prepared in state $\rho$. As we cannot measure precisely the quantum
state in a single experiment, we perform measurements on $n$ independent
identically prepared quantum systems. The measurement
carried out on each of the $n$ systems in state $\rho$ is done by QHT as
described in Section~\ref{QHT}. In the ideal setting, the results of such
experiments would be $n$  independent identically
distributed random variables $(X_{1}, \Phi_{1}),\dots ,(X_{n}, \Phi_{n})$  with
values in $\mathbb{R}\times [0,\pi]$ and distribution $P_{\rho}$ having density 
with respect to $\lambda$, ($\lambda$ being the Lebesgue measure on
$\mathbb{R}\times [0,\pi]$) equal to
      \begin{equation}
      \label{eq:prhophi}
     p_{\rho}(x,\phi)=\frac 1\pi   p_{\rho}(x|\phi) =\frac 1\pi
\mathcal{R}[W_{\rho}] (x,\phi),
       \end{equation}
where $\mathcal{R}$ is the Radon transform defined in equation
(\ref{eq.Radon.transform}).
As underlined in Section~\ref{QHT},  we do not observe $(X_{\ell},
\Phi_{\ell})_{\ell=1,\ldots n}$ but the noisy version $(Y_{\ell},
\Phi_{\ell})_{\ell=1,\ldots n}$  where
      \begin{equation}
      \label{noisy.data}
       Y_\ell:=\sqrt{\eta}\,X_\ell+\sqrt{(1-\eta)/2}\,\xi_\ell.
      \end{equation}
Here $\xi_\ell$'s are independent standard Gaussian random variables, independent
of all $(X_{\ell}, \Phi_{\ell})$, $\ell=1,\ldots, n$. The detection efficiency
$\eta\in]0.5,1]$ is a known parameter and $1-\eta$  represents the proportion of
photons which are not detected due to various losses in the measurement process.
\\
Let us denote by $p_{\rho}^{\eta}(y,\phi)$ the density of $(Y_{\ell},
\Phi_{\ell})$. Then,  for $\Phi=\phi$, the conditional density $p_{\rho}^{\eta}(\cdot|\phi)$
is the convolution of the density 
$\frac{1}{\sqrt{\eta}}p_\rho(\frac{\cdot}{\sqrt{\eta}}| \phi)$ of $\sqrt{\eta}X$
with $N^\eta$ the density of a centered Gaussian distribution having variance $(1
-\eta)/2$, that is 
       \begin{eqnarray}
       \label{densitbrui}
        p^\eta_\rho(y|\phi) &=& \left(\frac{1}{\sqrt\eta}
p_\rho\left(\frac{\cdot}{\sqrt\eta}|\phi\right) \ast N^\eta\right)(y)\nonumber\\
       & =& \int_{-\infty}^\infty \frac{1}{\sqrt{\eta}}
p_\rho\left(\frac{y-x}{\sqrt{\eta}}|\phi
        \right) N^\eta (x)dx.
        \end{eqnarray}
For $\Phi=\phi$, a useful equation in the Fourier domain, deduced by 
the previous relation (\ref{densitbrui}) is
       \begin{equation}
       \label{fourierproun}
       \mathcal{F}_1[\sqrt{\eta}p^\eta_\rho(\cdot\sqrt{\eta}|\phi)](t)
        = \mathcal{F}_1[p_\rho(\cdot|\phi)](t)
        \widetilde{N}^\gamma(t),
       \end{equation}
where $\mathcal{F}_{1} $ denotes the Fourier transform with respect to the first
variable and  $\widetilde{N}^\eta(t)=e^{-\frac{1-\eta}{4\eta}t^2}$ is the Fourier
transform of $N^\gamma(x)$, the density of a centered Gaussian density having
variance $(1 -\eta)/2\eta=\gamma$.\\ 

\noindent  In order  to estimate the elements of the density matrix defined in (\ref{rhojk})
from the data  $(Y_{\ell}, \Phi_{\ell})_{\ell=1,\ldots n}$, we define a
realistic class of quantum states $\mathcal{R}(C,B,r)$. For $C\geq 1$, $B>0$ and
$0 < r \leq 2$, the class $\mathcal{R}(C,B,r)$ is defined as follow
          \begin{equation}\label{eq.classcoeff}
         \mathcal{R}(C,B,r) \deq \{\rho {\rm \ quantum \ state} :
        |\rho_{m,n} |\leq C \exp(-B (m+n)^{r/2})\}.
        \end{equation}
Note that the class  $\mathcal{R}(C,B,r)$ has been translated in terms of Wigner
functions in \cite{ABM}, where it has been proved that the fast decay of the
elements of the density matrix implies both rapid decay of the Wigner function
and of its Fourier transform.\\

\noindent However, on the contrary to previous works, we do not assume here that the
constants $r$, $B$ and $C$ are known. From now on we denote by $\langle \cdot,
\cdot \rangle$ and $\| \cdot \|$ the usual Euclidean scalar product and norm.

\subsection{Outline of the results}
\label{Outline.results}


\noindent This paper deals with the problem of adaptive estimation of density matrix $\rho$ in QHT when taking into account the detection losses occurring in the measurement, leading to an additional Gaussian noise in the measurement data.  In order to compute  the performance of our procedure in $\mathbb{L}_2$ risk, we defined in previous section a realistic class of quantum states $ \mathcal{R}(C,B,r)$ in which the elements of the density matrix decrease rapidly. From the physical point of view, all the states which have been produced in the laboratory up to date belong to such a class, and a more detailed argument can be found in \cite{Butucea&Guta&Artiles}, as to why this assumption is realistic and in  \cite{ABM} as how to translate this class in terms of associated Wigner functions.\\

\noindent The problem of reconstructing the quantum state of a light beam has been extensively studied in quantum statistics and physical literature. Methods for reconstructing a quantum state are based on the estimation of either the density matrix $\rho$ or the Wigner function $W_\rho$.\\
\noindent The estimation of the density matrix from averages of data has been considered in the framework of ideal detection ($\eta=1$) in \cite{DAriano.0, DAriano.2, DAriano.3, Artiles&Gill&Guta}. Max-likelihood methods have been studied in \cite{BDPS,Artiles&Gill&Guta,DMS,Guta} and procedure using adaptive
tomographic kernels to minimize the variance has been proposed in \cite{DP}.  In a more general case of an efficiency parameter $\eta$ belonging to the interval $]1/2, 1]$,  the estimation of the density matrix of a quantum state of light has been discussed in \cite{DAriano.1,DMS,DAriano.5} and considered in \cite{Richter} via the pattern functions for the diagonal elements.  The problem of goodness-of-fit testing in quantum statistics has been
considered in \cite{MezianiTest}. In this noisy setting, the latter paper derived a testing procedure from a projection-type estimator where the projection is done in $L_2$ distance on some suitably chosen pattern functions. \\
\noindent For the problem of pointwise estimation of the Wigner function, we mention the work \cite{Guta&Artiles} in the case of ideal detection, that corresponds to $\eta=1$, where a kernel estimator is given and its sharp minimax optimality over a class of Wigner functions characterised by their smoothness is established. 
The same problem in the noisy setting  $\eta\in]1/2, 1]$ was treated in \cite{Butucea&Guta&Artiles}, where the minimax rates were obtained. The estimation of a quadratic functional of the Wigner function, as an estimator of the purity, was explored in \cite{Meziani}. \\
\noindent Recently, the more general case $\eta\in]0, 1]$ was investigated in \cite{ABM}. The authors provided rates of convergence in $ L_2$ loss for both an estimator of the Wigner function and an estimator of the density matrix.   Interestingly, the rates are polynomial in the case $r=2$, whereas they are intermediate for $r\in]0,2[$, where intermediate means that they are slower than any power of $n$ but faster than any power of $\log n$. However, the physicists argue, that their machines actually have high detection efficiency, around $0.9$. So we do not deal in this paper with values of $\eta$ smaller than $1/2$. It is to be noted that the estimator proposed in \cite{ABM} depends on the knowledge of $B$ and $r$. This is a serious limitation since in practice, one will face situations where one wants to reconstruct a density matrix without assuming the knowledge of $B$ and $r$. This is known in statistics as "adaptive estimation". In the present work, we tackle the problem of adaptive estimation over 
the classes of quantum states $\{\mathcal{R}(C,B,r)\}$. Our estimator is actually a soft-thresholded version of the estimator in \cite{ABM} which allows us to reach adaptation.\\
 Coefficients thresholding is now a classical tool in statistics.
It was introduced in a series of papers \cite{DJ94,DJ95,DJKP}
in the context of function estimation via wavelets coefficients. We refer to \cite{Tsybakov1998} for a comprehensive introduction to thresholding and waveltes. 
These methods were extended to inverse problems \cite{Dinv,Dinv2,kola,cava1}, see \cite{cava2} for an introduction and a survey of the most recent results.\\

\noindent The remainder of the article is organized as follows. In Section \ref{sec.dens.mat}, we present our adaptive thresholding procedure and state our main theoretical results. In particular, we establish upper bounds on the $L_2$ risk of our procedure and its achieves the convergence rates over a broad family of set $\mathcal{R}(C,B,r)$ which have been obtained in \cite{ABM}. These bounds are nonasymptotic and hold true with large probability. The theoretical investigation is complemented by numerical experiments reported in Section~\ref{sec.simulations}. The proofs of the main results are defered to the Appendix.

\section{Density matrix estimation}
\label{sec.dens.mat}

\noindent We assume now $n$ independent identically distributed random pairs $(Y_i,\Phi_i)_{i=1,\ldots,n}$ are observed, where $\Phi_1$ is uniformly distributed in $[0,\pi]$ and the conditional density of $Y_1$ given $\Phi_1$ is $p^\eta_\rho$, cf~(\ref{densitbrui}). The goal is to estimate the density matrix $[\rho_{j,k}]_{j,k}$ defined by  (\ref{rhojk}) and to investigate the convergence rate of the proposed estimator. To achieve this goal, we follow the framework of \cite{ABM} by assuming that the quantum state $\rho$ is in some class $\mathcal{R}(C,B,r)$ defined in (\ref{eq.classcoeff}). The notable difference of the present setting is that the precise knowledge of $C$, $B$ and $r$ are not required by our estimating procedure.

\subsection{Adapted pattern functions}
\label{Pattern.functions}

In order to reconstruct the entries of the density matrix from the noisy
observations $(Y_{\ell}, \Phi_{\ell})$ by a projection type estimator on the
\textit{pattern} functions, we have to adapt the pattern functions as follows. 
From now on, we shall use the notation $\gamma=\gamma(\eta) \deq
\frac{1-\eta}{4\eta}$.  We denote by $f^\eta_{k,j}$ the function which has
the following Fourier transform:
     \begin{equation}
      \label{eq:patterneta}
      \tilde{f}^\eta_{k,j}(t) \deq \tilde{f}_{k,j}(t) e^{\gamma t^2},
      \end{equation}
where $ \tilde{f}_{k,j}$ are the pattern functions defined in equation
(\ref{Rpattern}).

\subsection{Estimation procedure}
\label{Estimation.procedure}

The estimation procedure we introduce in this section will depend on one tuning parameter $N \deq N(n)$, the precise value of which will be given later. We define the set of indices
$J(N)\subset\mathbb{N}^2$ by
 \begin{equation}
 \label{eqJ}
J(N)\deq\{(j,k)\in\mathbb{N}^2,0\leq j+k \leq N-1\}.
\end{equation}
\noindent We first define an initial estimator $ \hat{\rho}^\eta$ of  $\rho$ by setting
      \begin{eqnarray}
      \label{estrhojk}
       \hat{\rho}^\eta_{j,k} \deq \left\{
       \begin{array}{ccc}
       \frac{1}{n}
\sum_{\ell=1}^{n}G_{j,k}\paren{\frac{Y_\ell}{\sqrt\eta},\Phi_\ell} &
\forall(j,k)\in J (N),\\
                                                                                
                                               0      & \mathrm{  otherwise, }
       \end{array}\right.
      \end{eqnarray}
where  $(G_{j,k})_{j,k}$ are constructed using the pattern functions
in~(\ref{eq:patterneta}) and 
        \begin{equation}
        \label{Gjk}
         G_{j,k}(x,\phi) \deq f^\eta_{j,k}(x) e^{-i(j-k)\phi}.
         \end{equation}
Note that this procedure introduced by  \cite{ABM} estimates the matrix coefficients by replacing the theoritical by its empirical conterpart.
To define our final procedure of estimation, let us introduce some notation.
From now, we denote by $\|.\|_{\infty}$ the supremum norm for functions, i.e. for any $f$,
$$\|f\|_{\infty}=\sup_{x\in\mathbb{R}}|f(x)|.$$
\noindent  Let $\varepsilon\in(0,1)$ be a prescribed tolerance level. The final estimation procedure applies the soft-thresholding operator to the initial one :
        \begin{eqnarray}
        \label{rho-seuil}
        \tilde{\rho}^{\eta}_{j,k} =  \frac{\hat{\rho}_{j,k}^{\eta}}{|\hat{\rho}_{j,k}^{\eta}|}
         \left(|\hat{\rho}_{j,k}^{\eta}|- t_{j,k} \right)_{+}, 
        \end{eqnarray}
        with the convention $0/0=0$, and
where the thresholds are defined as 
       \begin{equation}
      \label{tjk}
       t_{j,k} = 2\|f_{j,k}^{\eta}\|_{\infty}\sqrt{\frac{
\log\left(\frac{2N(N+1)}{\varepsilon}\right)}{n}}. 
      \end{equation}
Thus, our estimator of the density matrix is given by
$$\tilde{\rho}^{\eta}=[\tilde{\rho}^{\eta}_{j,k}]_{j,k} .$$

\subsection{Main results}
\label{Main.results}

To characterize the behaviour of the estimator $\tilde{\rho}^\eta$, we measure the quality of estimation in $\ell_{2}$-norm. For any density matrix $\nu=(\nu_{j,k})_{j,k\geq 0}$, we define the
$\ell_{2}$-norm of $\nu$ as
$$
 \|\nu \|_{2} = \sqrt{ \sum_{j,k\geq 0}| \nu_{j,k}|^{2} } .
 $$
We first state a risk bound that holds with large probability and will allow us to obtain the rates of convergence on the classes $\mathcal{R}(C,B,r)$.

 \begin{proposition}
           \label{thm.oracle}
           With probability at least $1-\varepsilon$, we have
           $$
          \left\|\tilde{\rho}^\eta - \rho \right\|^{2}_{2}  \leq
\inf_{I\subseteq J(N)} \left\{ 4 \sum_
            {(j,k)\in I} t_{j,k}^{2} + \sum_{(j,k)\notin I}|\rho_{j,k}|^{2}
\right\},
           $$
           where the set $J(N)$ is defined in (\ref{eqJ}).
\end{proposition}
The proof is given in the~\ref{sec.proofs.th}. Note that this result holds true for any value of the tuning parameter $N$. Choosing this parameter in a suitable manner leads to a rate of convergence that coincides with the one obtained in \cite{ABM} for a nonadaptive procedure. This result is stated in the following Theorem.

\begin{theorem}
           \label{coro1}
           Let us put $r_{0}\in(0,2)$, $B_{0}>0$ and let us choose
                    \begin{equation}
         \label{Ncorr1}
                      N = N(n) := \left\lfloor \left(\frac{\log(n)}{2B_{0}}\right)^{\frac{2}{r_{0}}}
                   \right\rfloor,
           \end{equation}
where $\left\lfloor x \right\rfloor$ denotes the integer part of $x$ such that $\left\lfloor x \right\rfloor\leq x<  \left\lfloor x \right\rfloor+1$. Let us assume that $\rho\in\mathcal{R}(C,B,r)$, for some
unknown $C\geq 1$, $B\geq B_{0}$, $r\in[r_0 , 2]$. Then, there are constants
$\mathcal{C}_{1},\mathcal{C}_ {2},\mathcal{C}_{3}>0$ such that with probability
at least $1-\varepsilon$, we have\\
\begin{itemize}
\item For  $\eta = 1$ and $r\in[r_0,2]$
            \begin{equation*}
                    \left\|\tilde{\rho}^\eta - \rho \right\|^{2}_{2}
\leq\mathcal{C}_1 n^{-1}\left(\log(n)\right)^{\frac{20}{3r}}\log
\left(\log(n)\varepsilon^{-1}\right).
              \end{equation*}    \noindent \\
\item   For $\eta\in(\frac{1}{2},1)$ and $ r=2$             
            \begin{equation*}
                    \left\|\tilde{\rho}^\eta - \rho \right\|^{2}_{2}\leq
\mathcal{C}_2   n^{-\frac{B}{4\gamma+B}}\left(\log(n)+\left(\log(n)\right)^{1/3}
\log\left( \log (n)\varepsilon^{-1}\right) \right). 
            \end{equation*} \noindent \\               
\item   For $\eta\in(\frac{1}{2},1)$ and  $r\in(r_0,2)$                
            \begin{equation*}
                    \left\|\tilde{\rho}^\eta - \rho \right\|^{2}_{2}\leq 
\mathcal{C}_3   e^{-2BM(n)^{r/2}} \left( \log(n)^{2-r/2}
+\log(n)^{1/3}\log\left( \log (n)\varepsilon^{-1}\right)  \right),
                   \end{equation*}\noindent \\
 where $ M(n) $ satisfies $8 \gamma M(n) + 2 B M(n)^{r/2} = \log(n)$.  In
particular, note that 
 $$ M(n)= \frac{1}{8\gamma} \log(n) - \frac{2B}{(8\gamma)^{1+r/2}} \log(n)^{r/2}
+ o(\log(n)^{r/2}).$$
 \end{itemize}
\end{theorem}
The proof is given in the~\ref{sec.proofs.co1}. Let us give some comments on this result highlighting its relations to previous work. First of all, note that the convergence rate is polynomial in the cases $(\eta,r)\in\{2\}\times[r_0,2]$ and  $(\eta,r)\in(1/2,1)\times\{2\}$. Furthermore, the rate is parametric, up to a logarithmic factor, in the first case. It is slower in the second case, but becomes closer to the parametric rate when $B$ is very large. The benefits of the adaptation are particularly striking in this case. Indeed, if, for example, the only available information is that $B\geq 1/3$ and $\eta=3/4$, then the estimator proposed in \cite{ABM} will converge at the rate $n^{1/2}\log n$, even if the true state $\rho$ belongs to the class $\mathcal{R}(C,B,2)$ with a very large constant $B\geq 1/3$. Contrarily to this, our estimator will converge at the rate $n^{-\frac{3B}{1+3B}}\log n$ which can be very close to the parametric rate $n^{-1}$ if $B$ is large.\\
\noindent One can also note that when $(\eta,r)\in(1/2,1)\times(r_0,2)$, the rate we get is slower than any power of $n^{-1}$, but faster that any power of $(\log n)^{-1}$. We will say that these rates are intermediate. They coincide, up to a $\log(\log(n)/\varepsilon)$ factor, with the rates obtained in \cite{ABM,Katia}. Another interesting feature of the previous result is that it provides a risk bound with high probability, whereas existing results are all concerned with bounding the expected risk.\\
Interestingly, the same procedure achieves the nearly parametric rate in the case of pure state as well.

\begin{theorem}
 \label{coro2}
              Under the same choice for $N$ in Theorem \ref{coro1},
              $$
               N = N(n) := \left\lfloor
\left(\frac{\log(n)}{2B_{0}}\right)^{\frac{2}{r_{0}}} \right\rfloor, 
              $$
            if $\rho$ is a pure state, i.e., if $\rho_{j_0,j_0} = 1$ for some $j_0$ and all the other $\rho_{j,k}$'s are $0$. Then we have, as soon as $N> \max(j_{0},2)$, with
probability at least $1-\varepsilon$, 
              $$
              \left\|\tilde{\rho}^\eta - \rho \right\|^{2}_{2}  < \frac{64}{n
r_0} \|f_{j_0,j_0}\|_{\infty}^2 \log\left(\frac{2 \log(n)} {B_0
\varepsilon}\right).
              $$
\end{theorem}
The proof is given in the~\ref{sec.proofs.co2}.

\newpage

\section{Experimental evaluation}
\label{sec.simulations}

\subsection{Examples considered in the experiments}
\label{example}

We present in Table~\ref{TTtab:1} examples of pure quantum states, which can be
created at this moment in laboratory and belong to the class $\mathcal{R}(C,B,r)$ with $r=2$.   Table~\ref{TTtab:1} gives also their density matrix coefficients $\rho_{j,k}$ and probability densities $p_\rho(x|\phi)$.\\
Among the pure states we consider the \textit{vacuum} state, which is the pure state of zero photons. Note that the \textit{vacuum} state would
provide a random variable of Gaussian probability density $p_\rho(x|\phi)$ via the ideal measurement of QHT (see Section~\ref{QHT}). That explains the Gaussian nature of the noise in the effective result of the QHT measurement. \\
 We consider also the \textit{single photon} state which is the pure state of one photon and the \textit{coherent-$q_0$} state, which characterizes the laser pulse with the number of photons  Poisson distributed with an average of $M$ photons. Remark that the well-known \textit{Schr\"{o}dinger cat} state is described by a linear
superposition of two \textit{coherent} vectors (see e.g. \cite{Our}).\\

\begin{table}[htbp!]
\begin{center}%
\caption{Examples of quantum states}
\label{TTtab:1}
\vspace{.2cm}
\begin{tabular}{|l|}
\hline
\textbf{\textit{Vacuum} state}  \\
$\bullet$ $\rho_{0,0}=1$ rest zero,\\
$\bullet$ $p_\rho(x|\phi)=e^{-x^2}/\sqrt{\pi}$. \\
\hline
\textbf{Single photon state}\\
$\bullet$ $\rho_{1,1}=1$ rest zero,\\
$\bullet$ $p_\rho(x|\phi)=2x^2e^{-x^2}/\sqrt{\pi}.$\\
\hline
\textbf{\textit{Coherent-$q_0$} state }$q_0\in\mathbb{R}$ \\
$\bullet$ $\rho_{j,k}=e^{-|q_0|^2}(q_0/\sqrt{2})^{j+k}/\sqrt{j!k!}$,\\
$\bullet$ $p_\rho(x|\phi)=\exp(-(x-q_0\cos(\phi))^2)/\sqrt{\pi}.$\\
\hline            
\textbf{\textit{Thermal state}} $\beta>0$\\
$\bullet$ $\rho_{j,k}=\delta^j_k(1-e^{-\beta})e^{-\beta k}$,\\
$\bullet$ $p_\rho(x|\phi)=\sqrt{\tanh(\beta/2)/\pi}\exp(-x^2\tanh(\beta/2)).$\\
\hline
\textbf{\textit{Schr\"{o}dinger cat}} $q_0>0$ \\
$\bullet$
$\rho_{j,k}=2(q_0/\sqrt{2})^{j+k}/\left(\sqrt{j!k!}
(\exp(q_0^2/2)+\exp(-q_0^2/2))\right)$, for $j$ and $k$ even, rest zero,\\
$\bullet$ $p_\rho(x|\phi)=\left(\exp(-(x-q_0\cos(\phi))^2)+\exp(-(x+q_0
\cos(\phi))^2)\right.$\\
$ \quad$ $ \quad$ $ \quad$ $ \quad$
$\left.+2\cos(2q_0x\sin(\phi))\exp(-x^2-q_0^2\cos^2(\phi))\right)
 /\left(2\sqrt{\pi}(1+\exp(-q_0^2))\right).$\\
 \hline
\end{tabular}
\end{center}
\end{table}

\subsection{Pattern functions $f^\eta_{j,k}$}
\label{TTsimu.pattern.function}

\begin{figure}[htbp!]
\caption{Examples of pattern functions $f_{j,k}$ (a) and adapted pattern functions $f_{j,k}^\eta$ (b).}
\label{fig:RMSE1}
        \centering
        \subfigure[]{\includegraphics[scale=0.52]{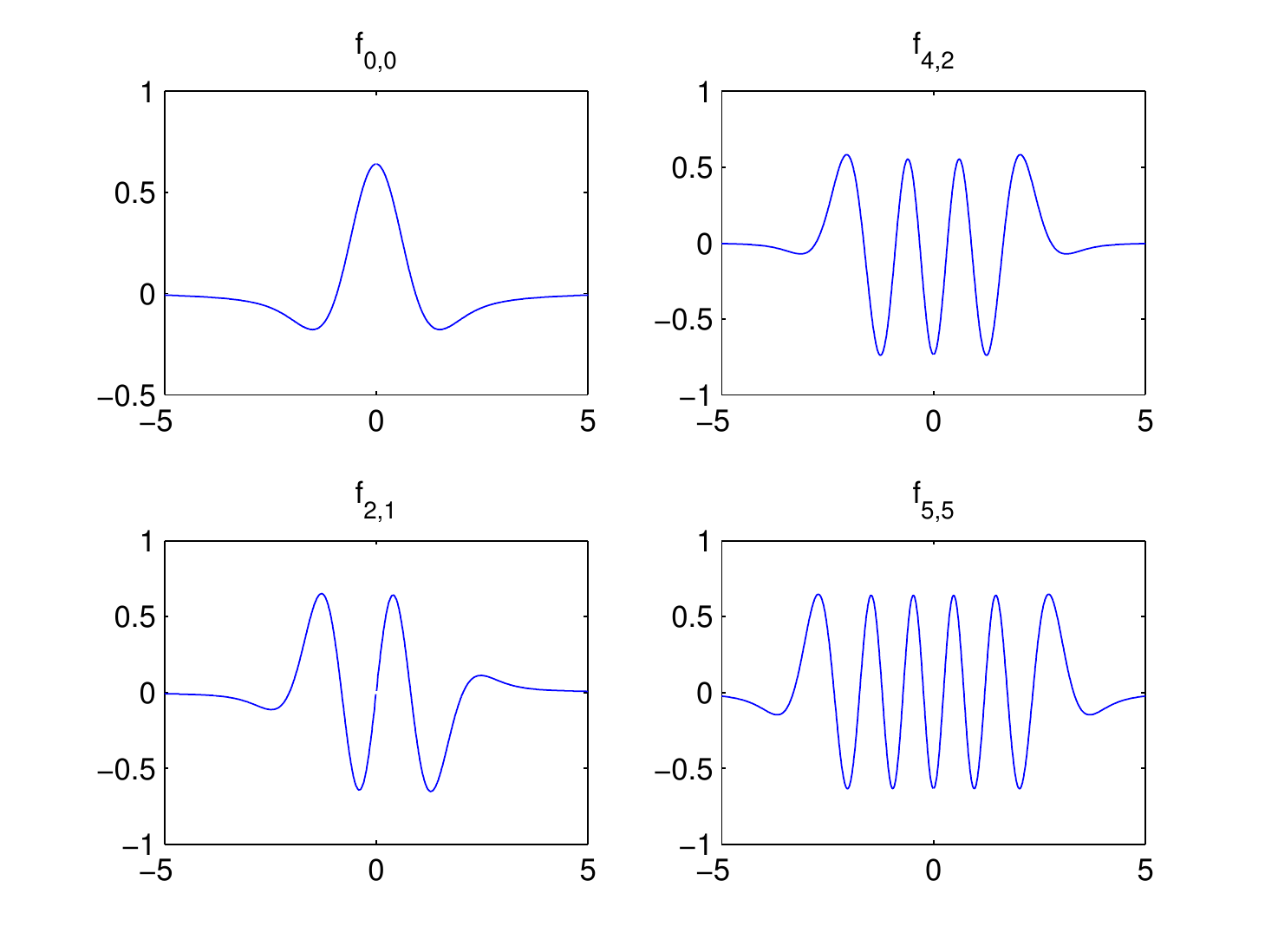}}
        \subfigure[]{\includegraphics[scale=0.52]{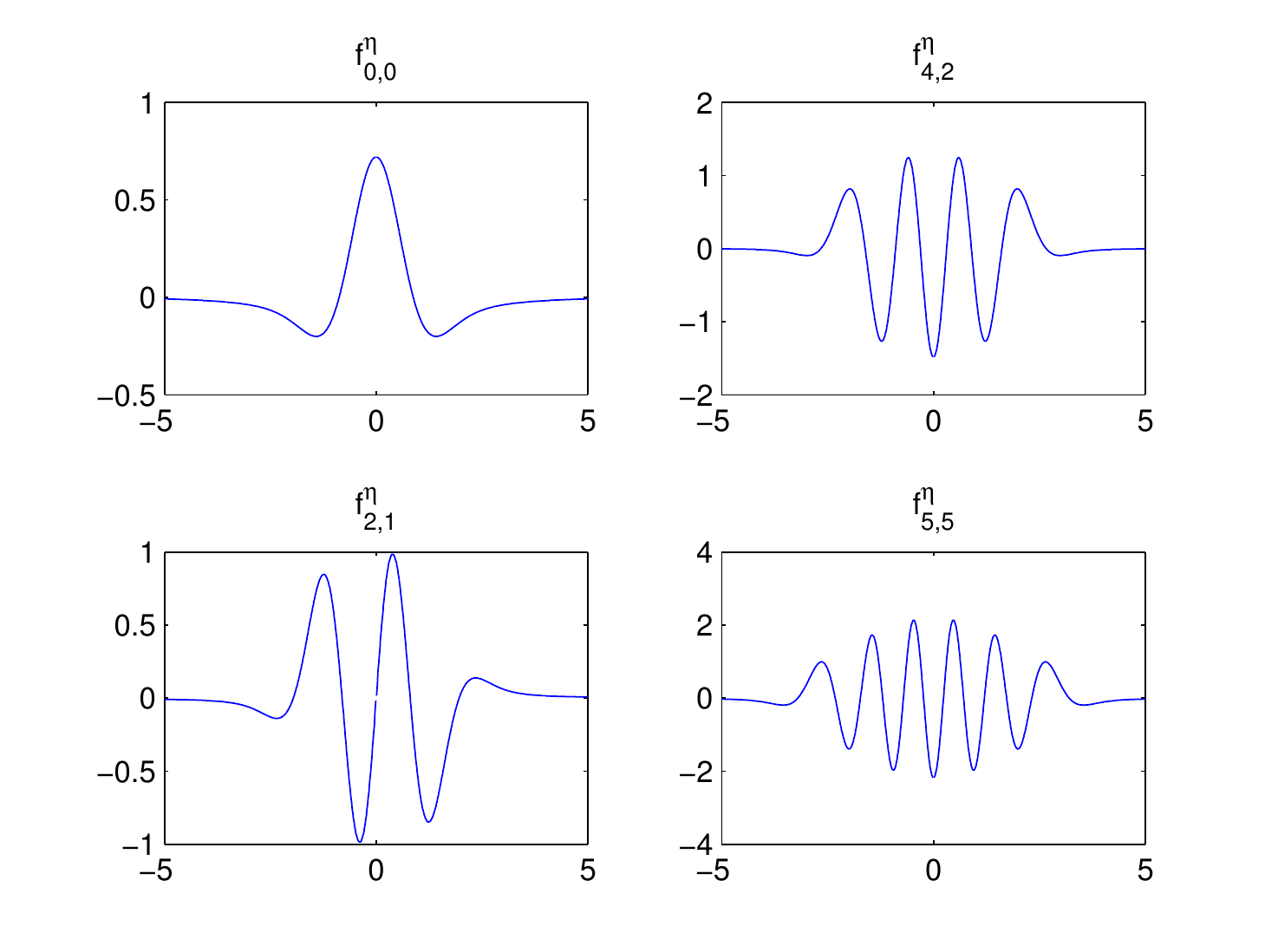} }
\end{figure}

Since there is no closed-form expression for the pattern functions $f^\eta_{j,k}$, we evaluate them numerically on a 1-D regular grid of $Q=4096$ points. We use expressions (\ref{Rpattern}) and (\ref{eq:patterneta}) to evaluate $\tilde f_{j,k}$ and $\tilde f_{j,k}^\eta$  on the 1-D frequency grid of $Q$ discretized $t$ points. The adapted pattern functions $f_{j,k}^\eta$ are computed on the 1-D spacial grid of $Q$ discretized $x$ points by applying to $\tilde f_{j,k}^\eta$ the inverse Fast Fourier Transform (FFT) in $O(Q \log(Q))$ operations. Some pattern and adapted functions are depicted in Figure~\ref{fig:RMSE1}.

\subsection{Implementation of our procedure}
\label{TTimplementation.Mn}

\begin{figure}[htbp!]
\caption{First row: The density matrix $\rho$ respectively of the coherent state, the chr\"{o}dinger cat state and Thermal state.
  Following rows: estimated $\tilde \rho^\eta$ of  previous states for $B_0=0.5$, $\eta=0.9$, $\varepsilon=1$ and $n$ respectively equal to 
  $10 \times 10^3$ (row \#2), $100 \times 10^3$ (row \#3), $500 \times 10^3$ (row \#4).}
\label{fig:examples}
        \centering
        \subfigure[Coherent $q_0=3$]{\includegraphics[scale=0.35]{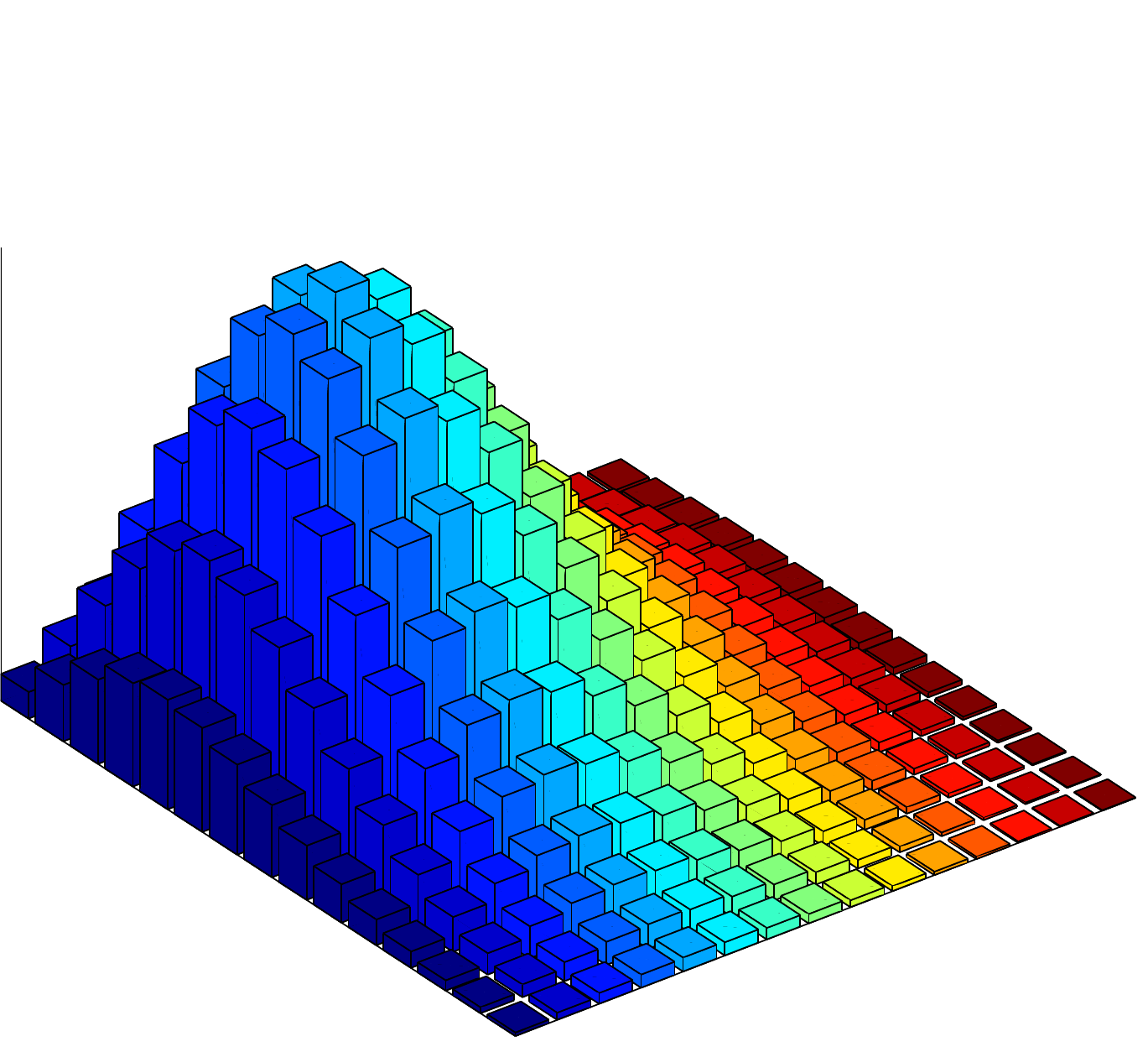}}\,\,
        \subfigure[Schr\"{o}dinger cat $q_0=3$ ]{\includegraphics[scale=0.35]{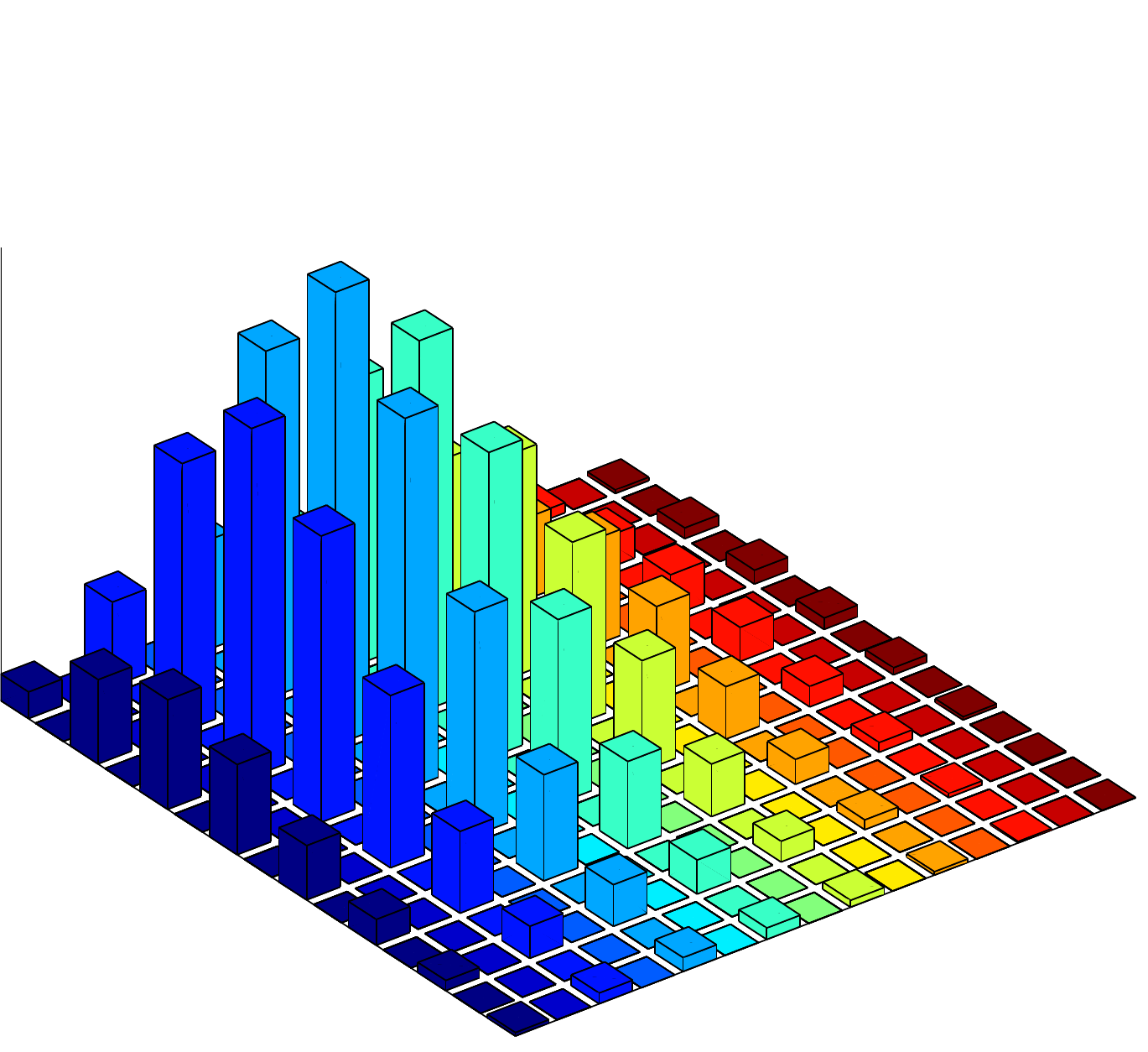} }\,\,
         \subfigure[Thermal $\beta=1/4$]{\includegraphics[scale=0.35]{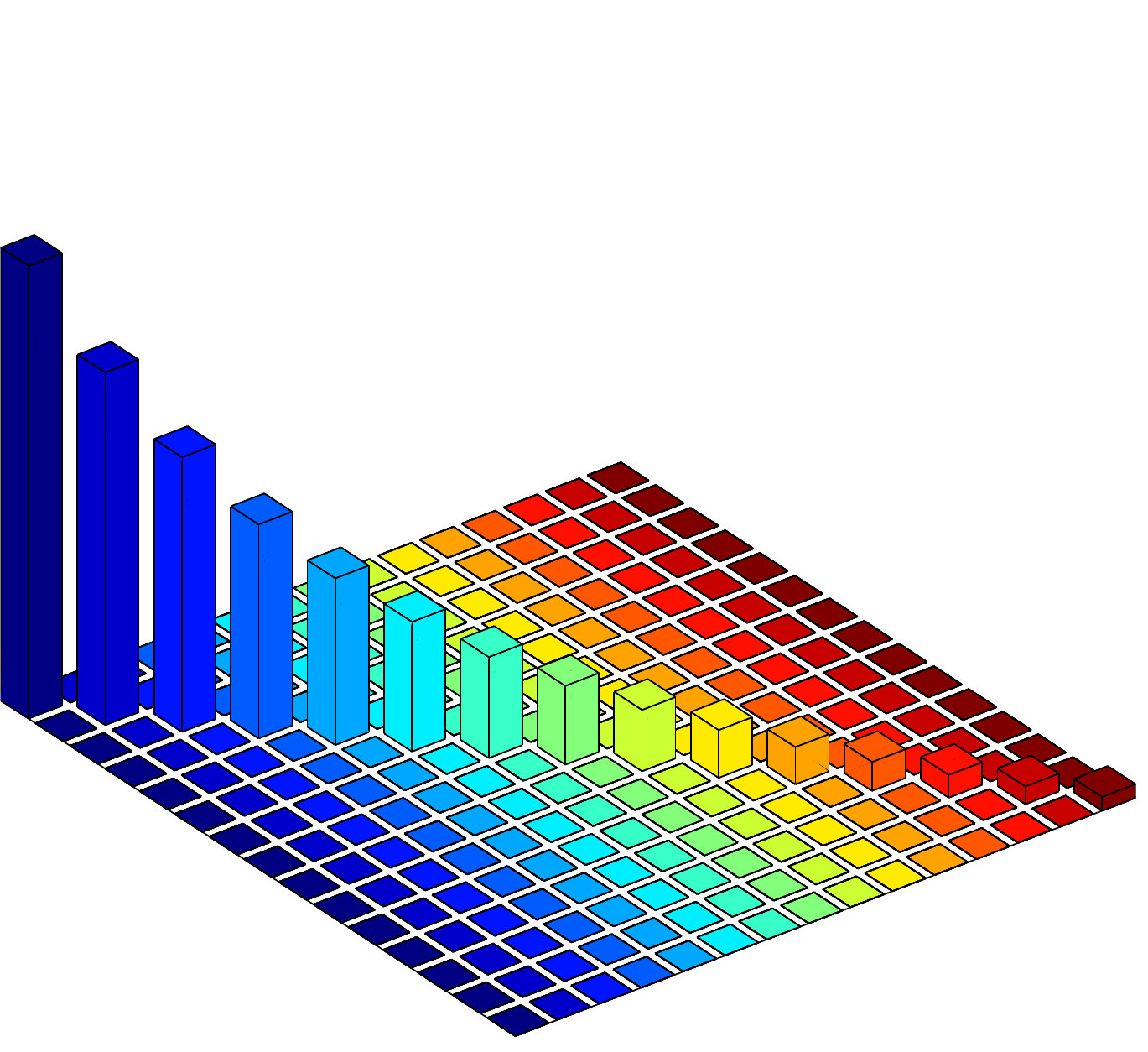} }
            \centering
        \subfigure{\includegraphics[scale=0.35]{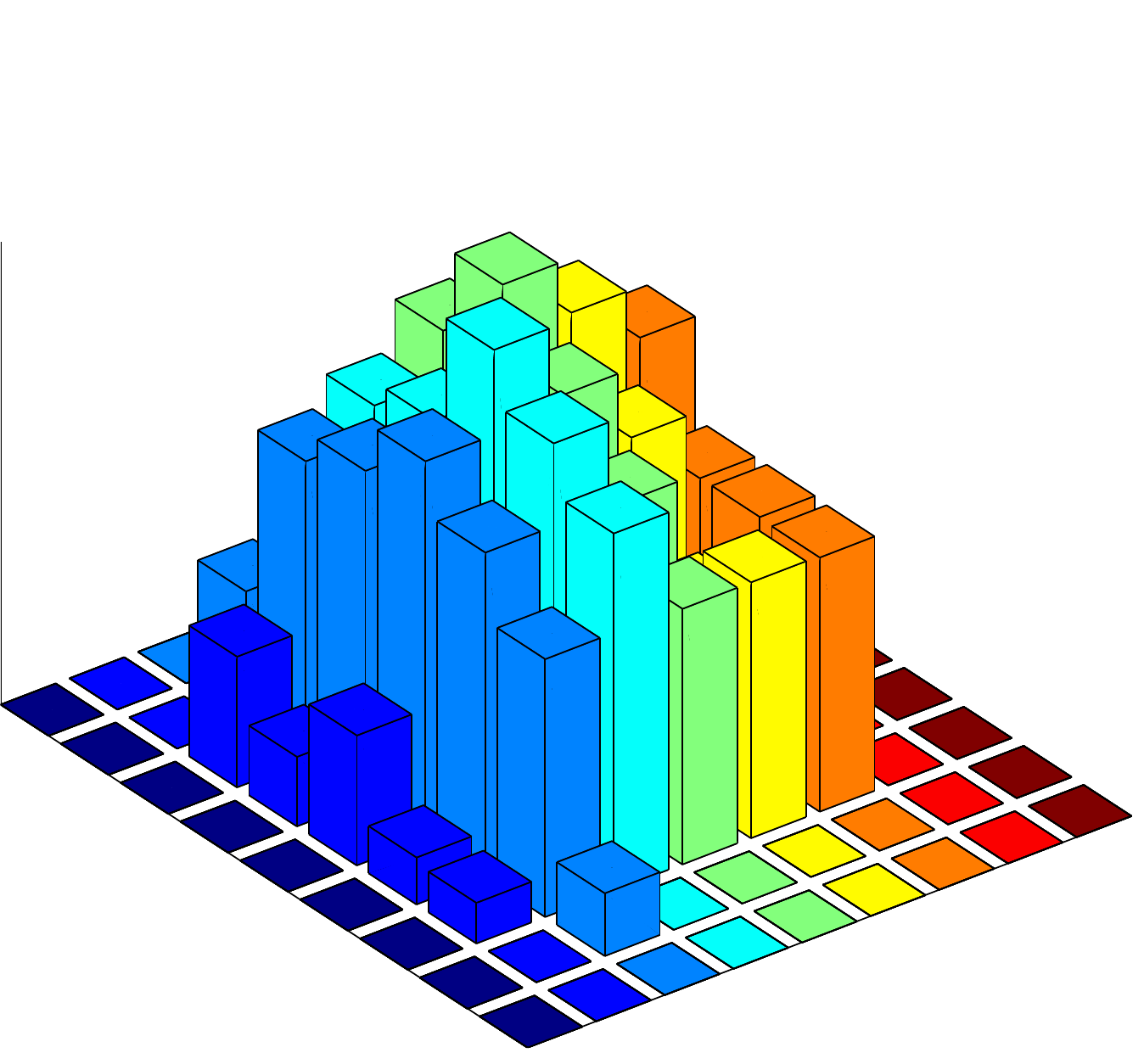}}\,\,
        \subfigure{\includegraphics[scale=0.35]{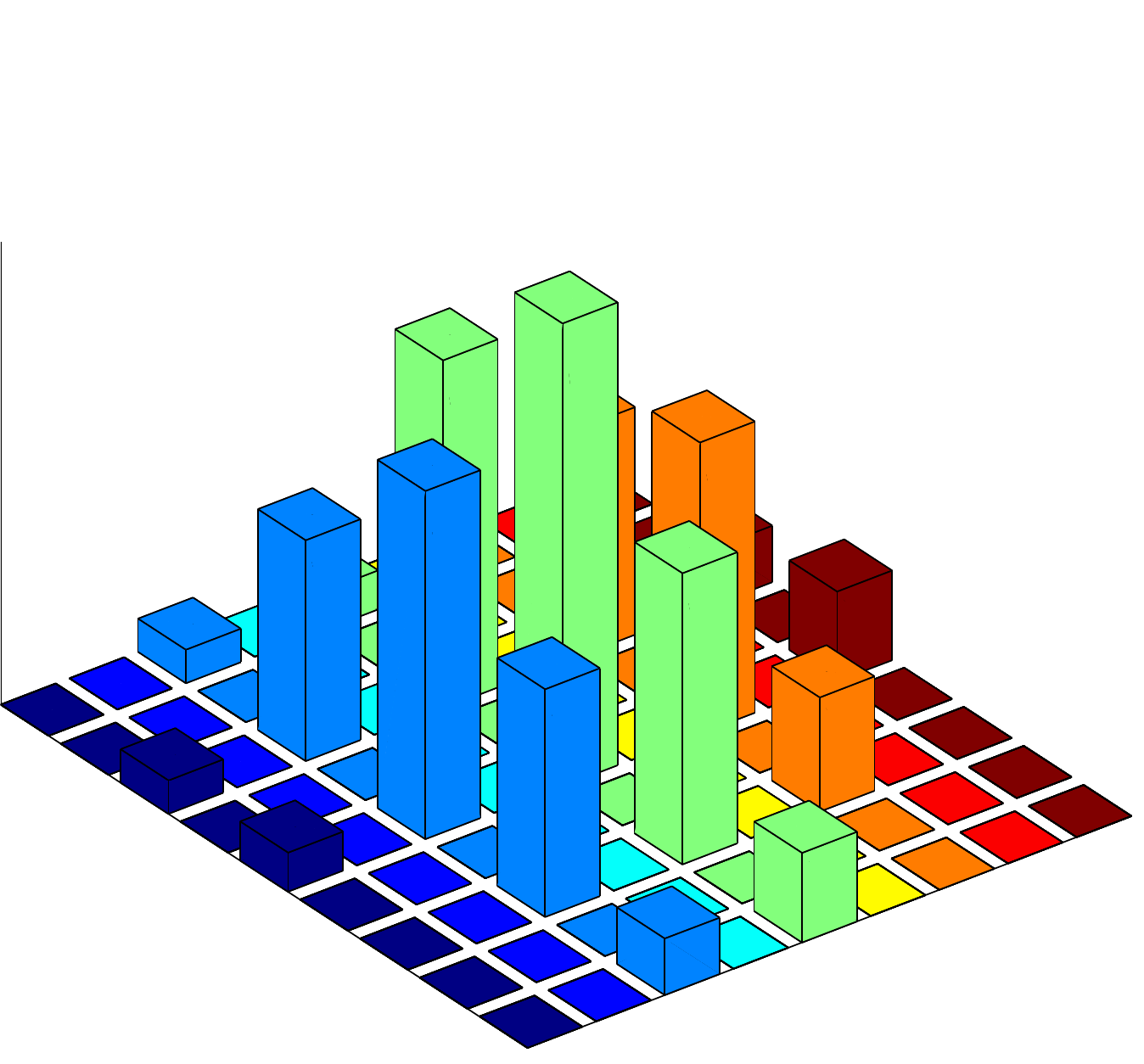} }\,\,
         \subfigure{\includegraphics[scale=0.35]{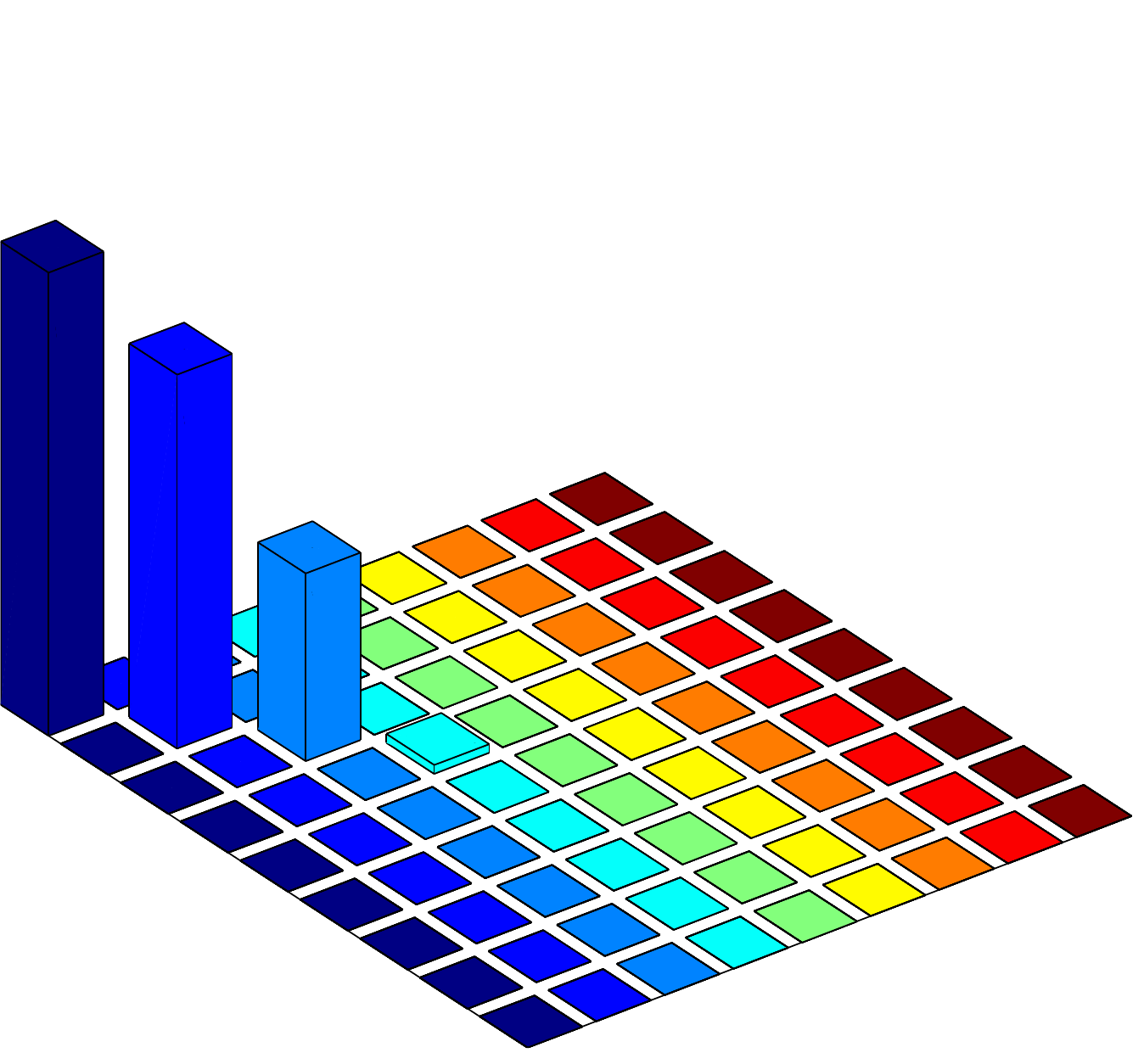} }     
            \centering
        \subfigure{\includegraphics[scale=0.35]{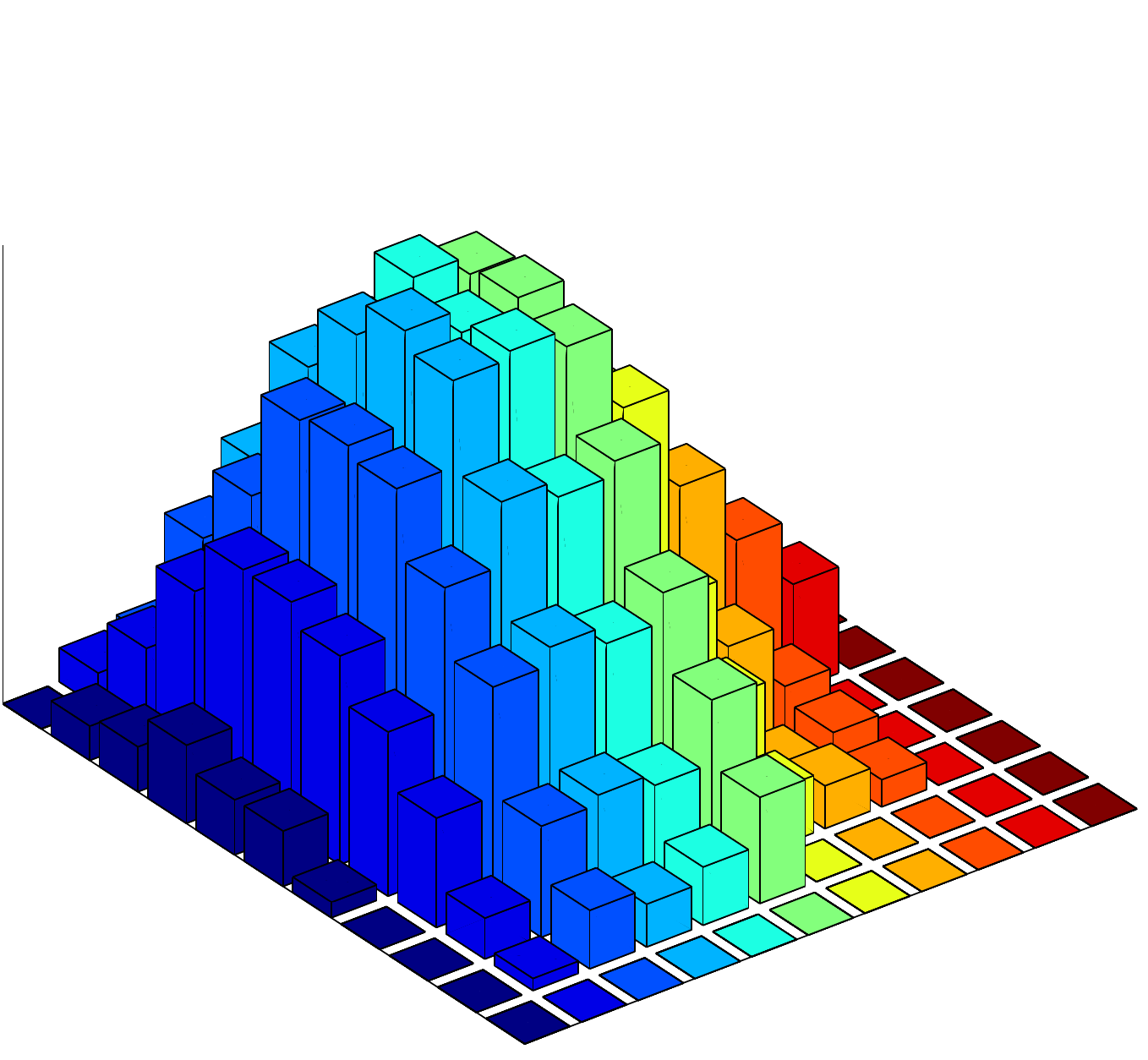}}\,\,
        \subfigure{\includegraphics[scale=0.35]{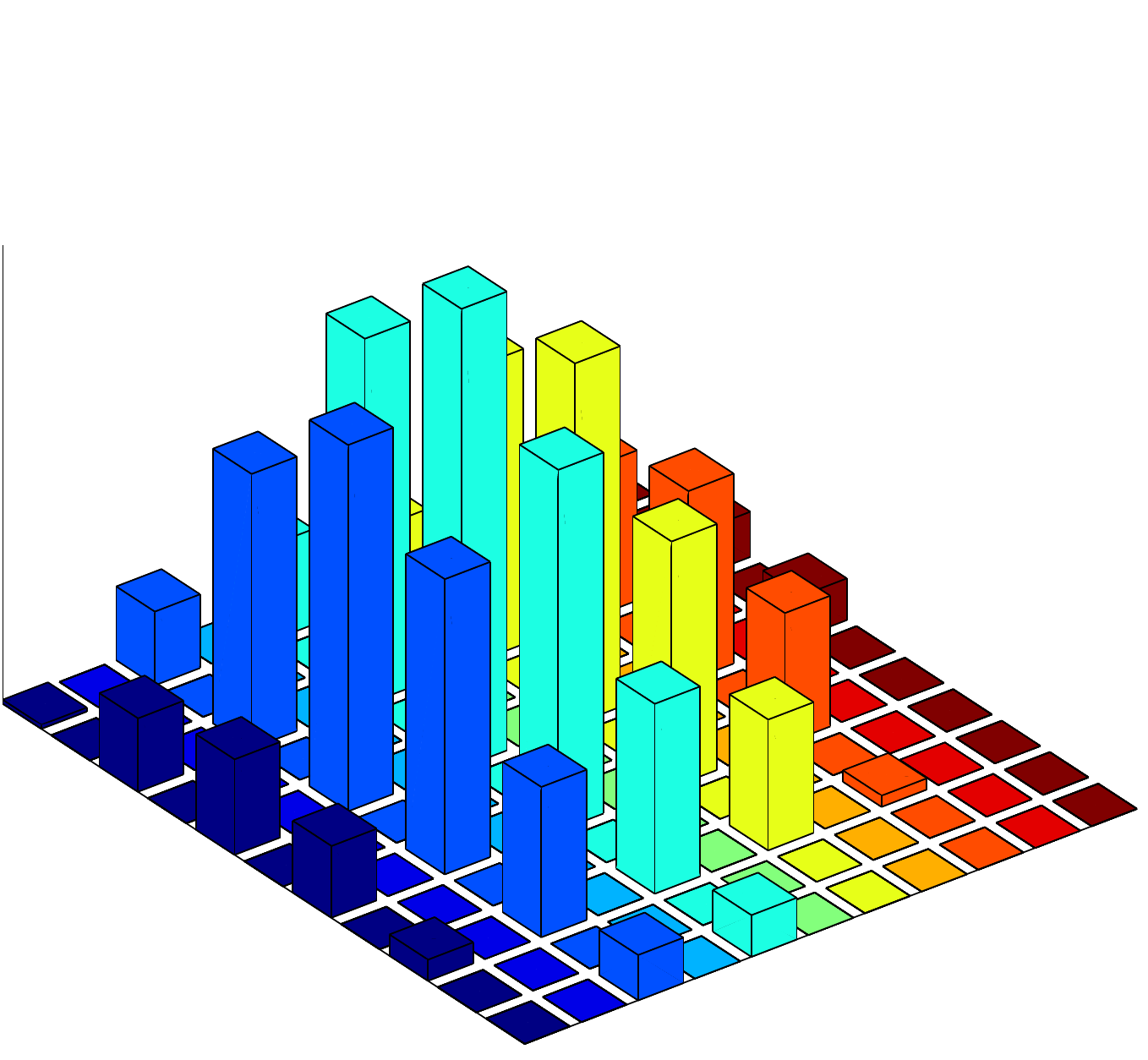} }\,\,
         \subfigure{\includegraphics[scale=0.35]{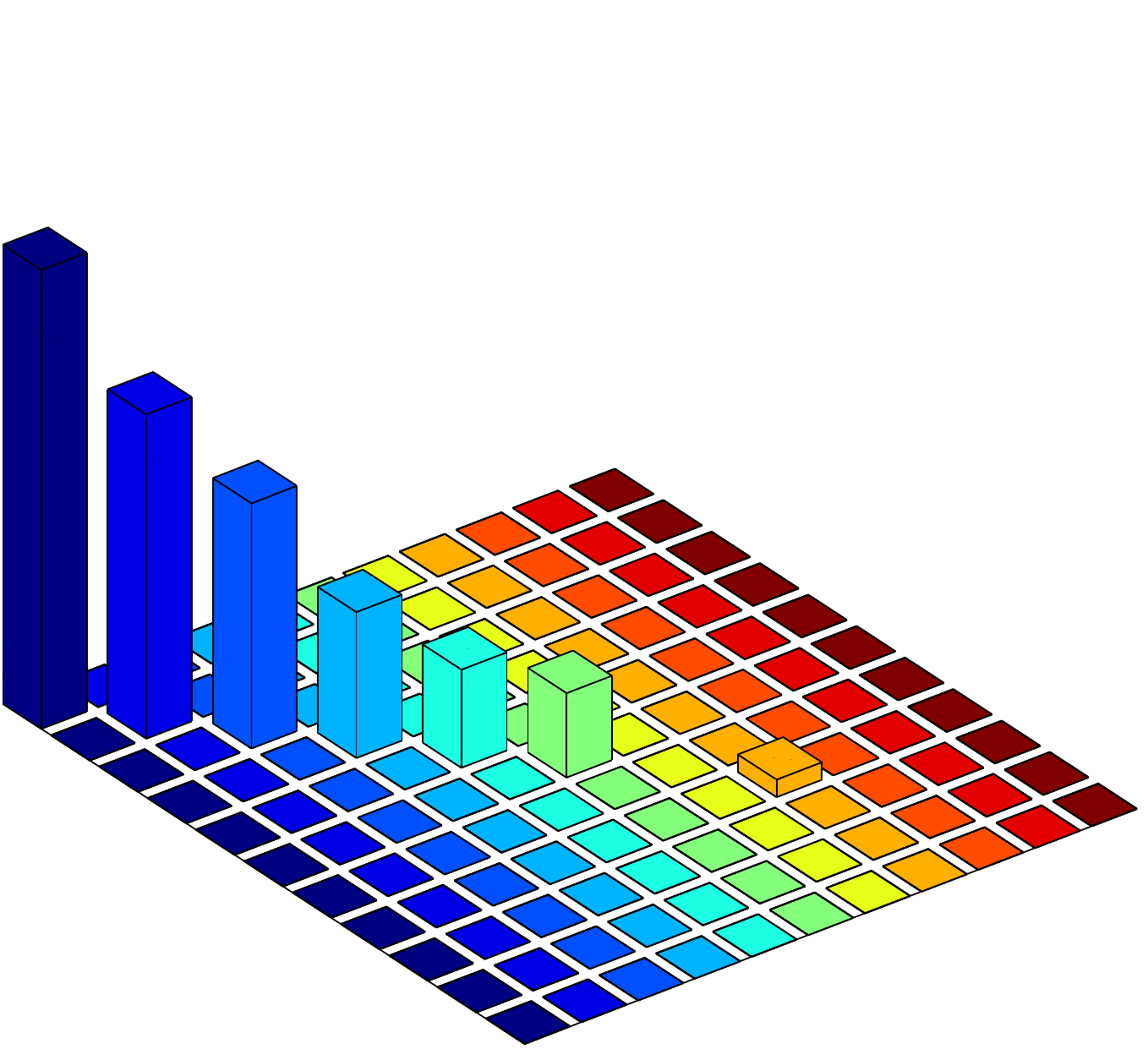} }   
            \centering
        \subfigure{\includegraphics[scale=0.35]{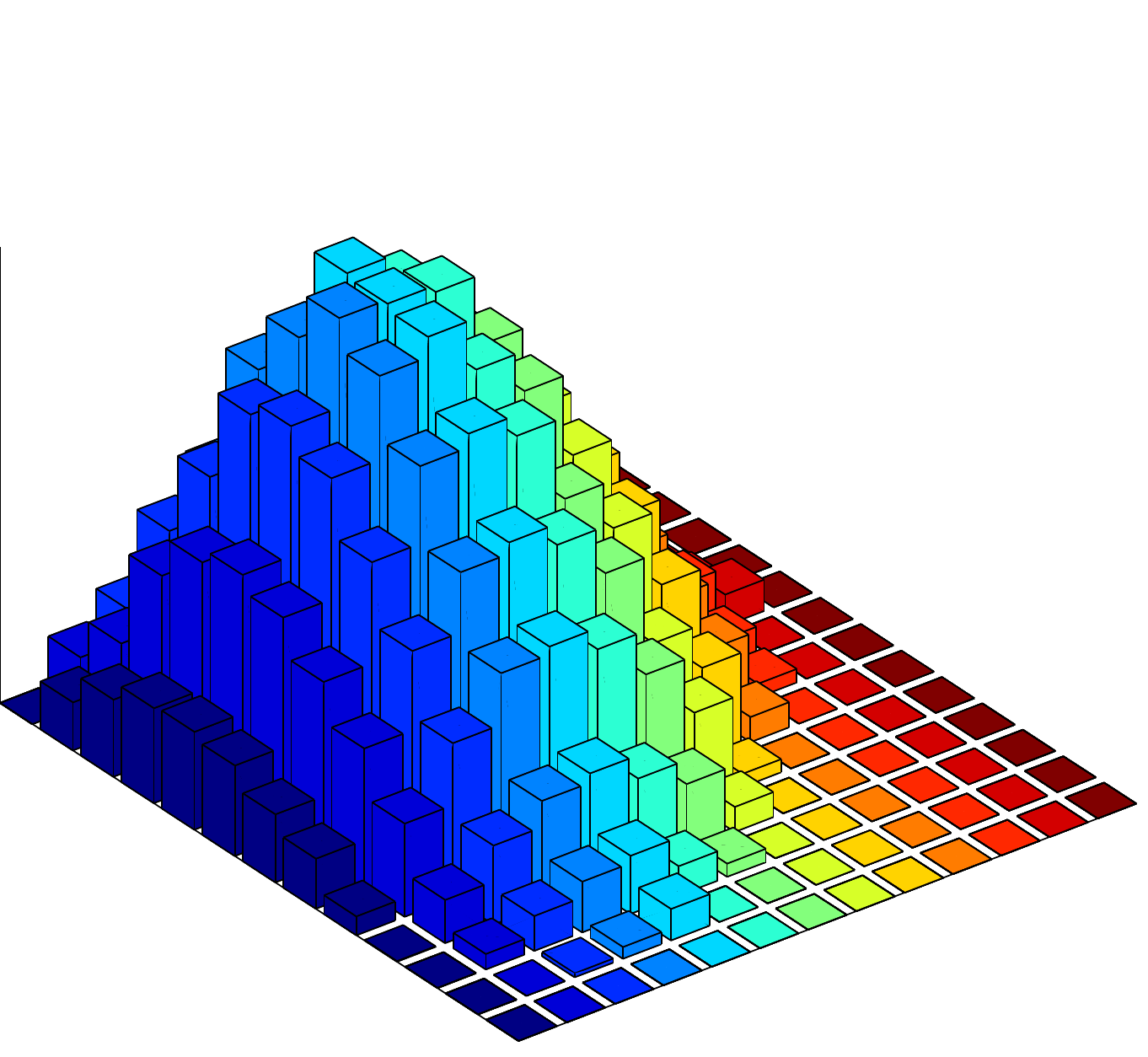}}\,\,
        \subfigure{\includegraphics[scale=0.35]{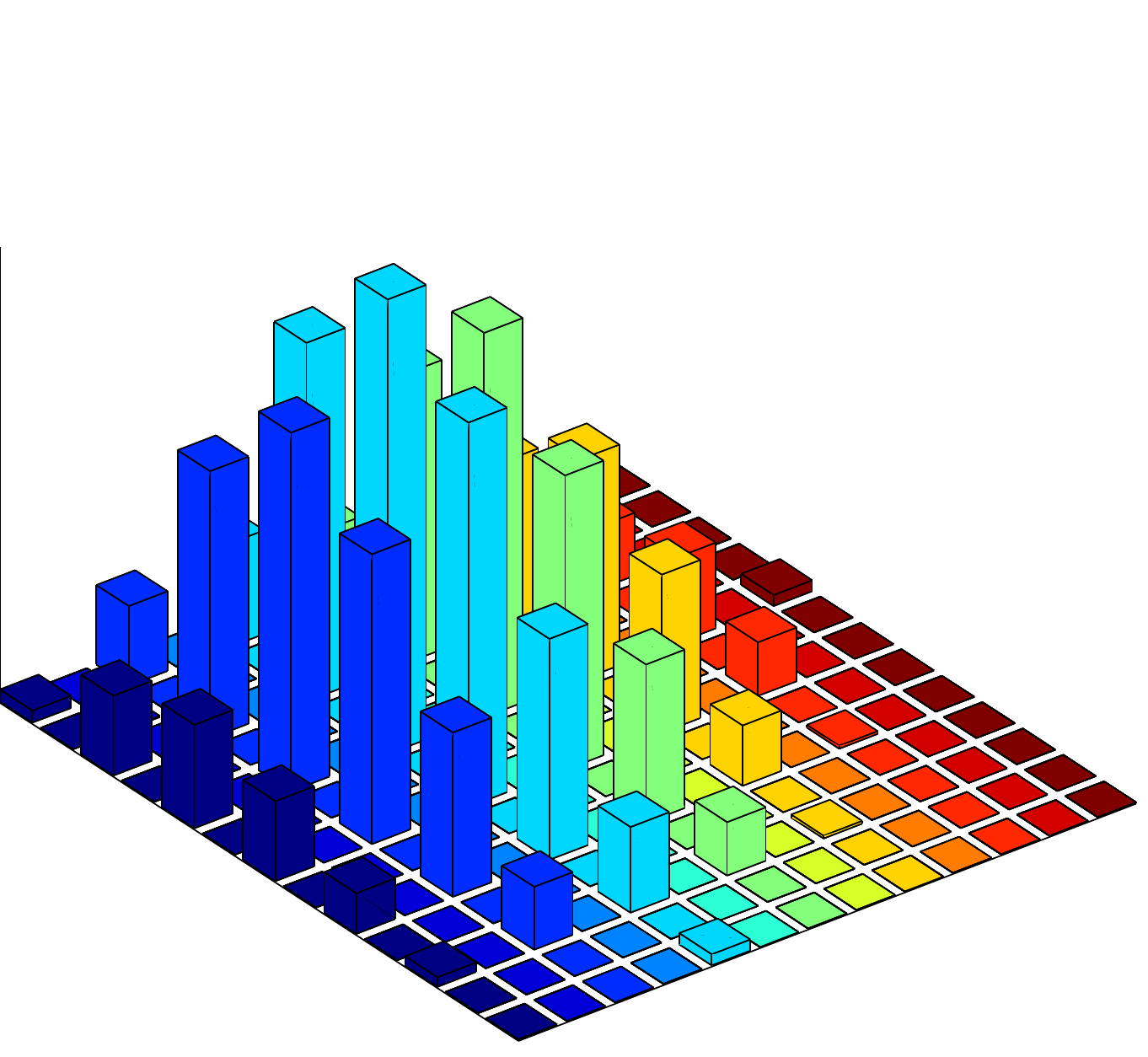} }\,\,
         \subfigure{\includegraphics[scale=0.35]{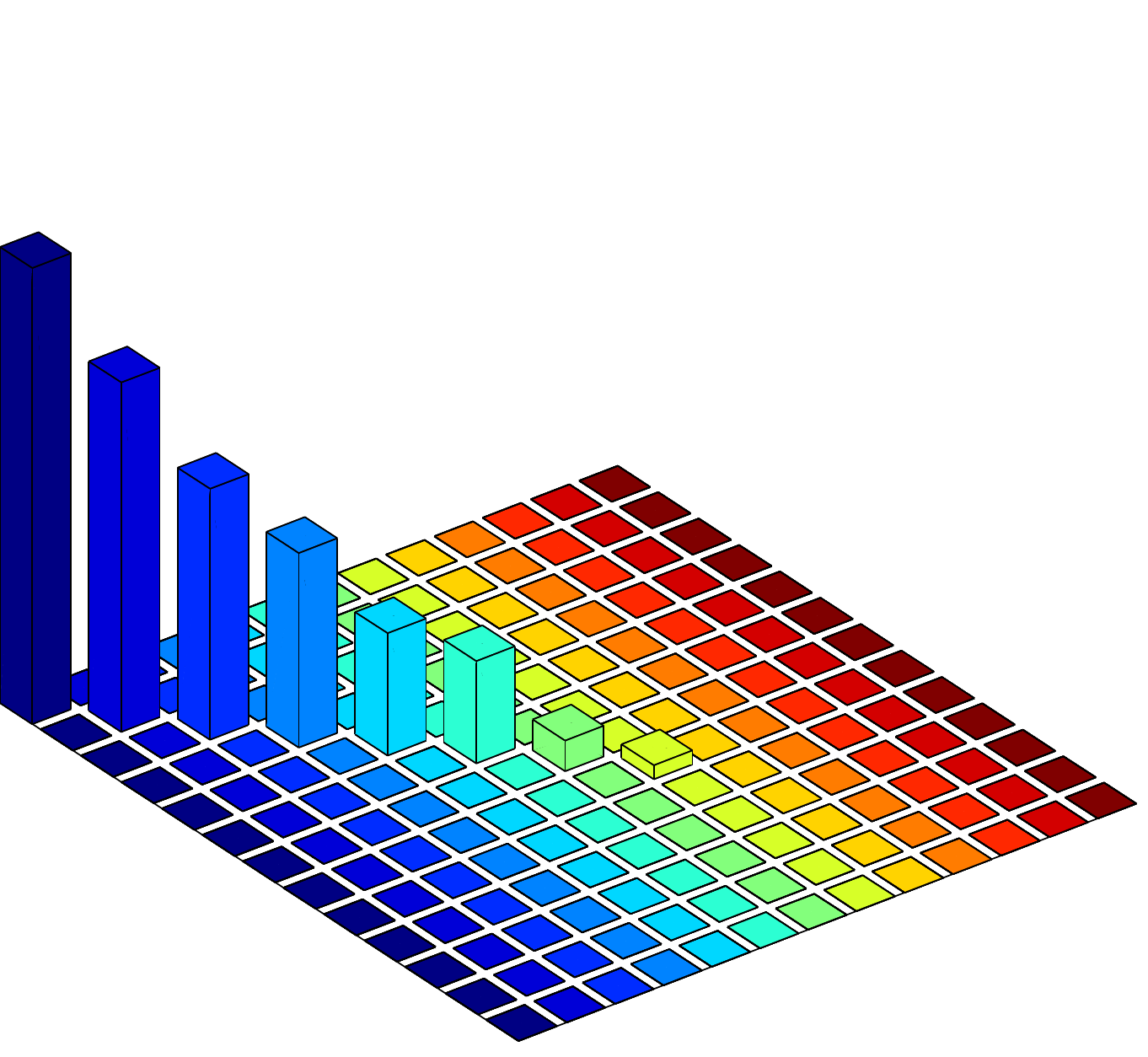} }   
       
\end{figure}
The deconvolved estimator $\hat{\rho}^\eta_{j,k}$ defined in (\ref{estrhojk}) is computed by evaluating 
$$ G_{j,k}(x,\phi) = f^\eta_{j,k}(x) e^{-i(j-k)\phi}$$
at point $x$ using a cubic spline interpolation of the values of $f^\eta_{j,k}$ on the discrete grid of $Q$ points.\\

\noindent In the following section, we assess the performance of the threshold estimator $\tilde{\rho}^{\eta}_{j,k}$. We perform this evaluation by creating noisy samples $Y_\ell$ as defined in (\ref{noisy.data}). The initial samples $X_\ell$ are drawn from the distribution $p_\rho(x | \phi)$ (see Table \ref{TTtab:1}) using the rejection method. The value of $N=N(n)$ is set following (\ref{Ncorr1}). We use $r_0=2$ and $B_0=1/2$ for all the numerical experiments. A toolbox that implements this procedure and reproduces all the figures of this article is available online\footnote{\textit{http://www.ceremade.dauphine.fr/~peyre/codes/}}. In Figure~\ref{fig:examples}, represents the density matrices $\rho$ and the estimated density matrices  $\tilde \rho^\eta$ of some quantum state.

\newpage
\subsection{Studies of the performance of our estimation procedure }
\label{TTperformanceOmega}

\begin{figure}[htbp!]
\caption{Evolution of $\mathbb{E}( \textrm{RMSE}(n) )$ as a function of $n$ for $\eta=0.9$, $\varepsilon=1$ and $N=30$.
	The blue shaded area represent the confidence interval 
	at $\pm 3$ times the standard deviation of RMSE$(n)$.}
\label{fig:RMSE}
        \centering
        \subfigure[Coherent $q_0=3$]{\includegraphics[scale=0.42]{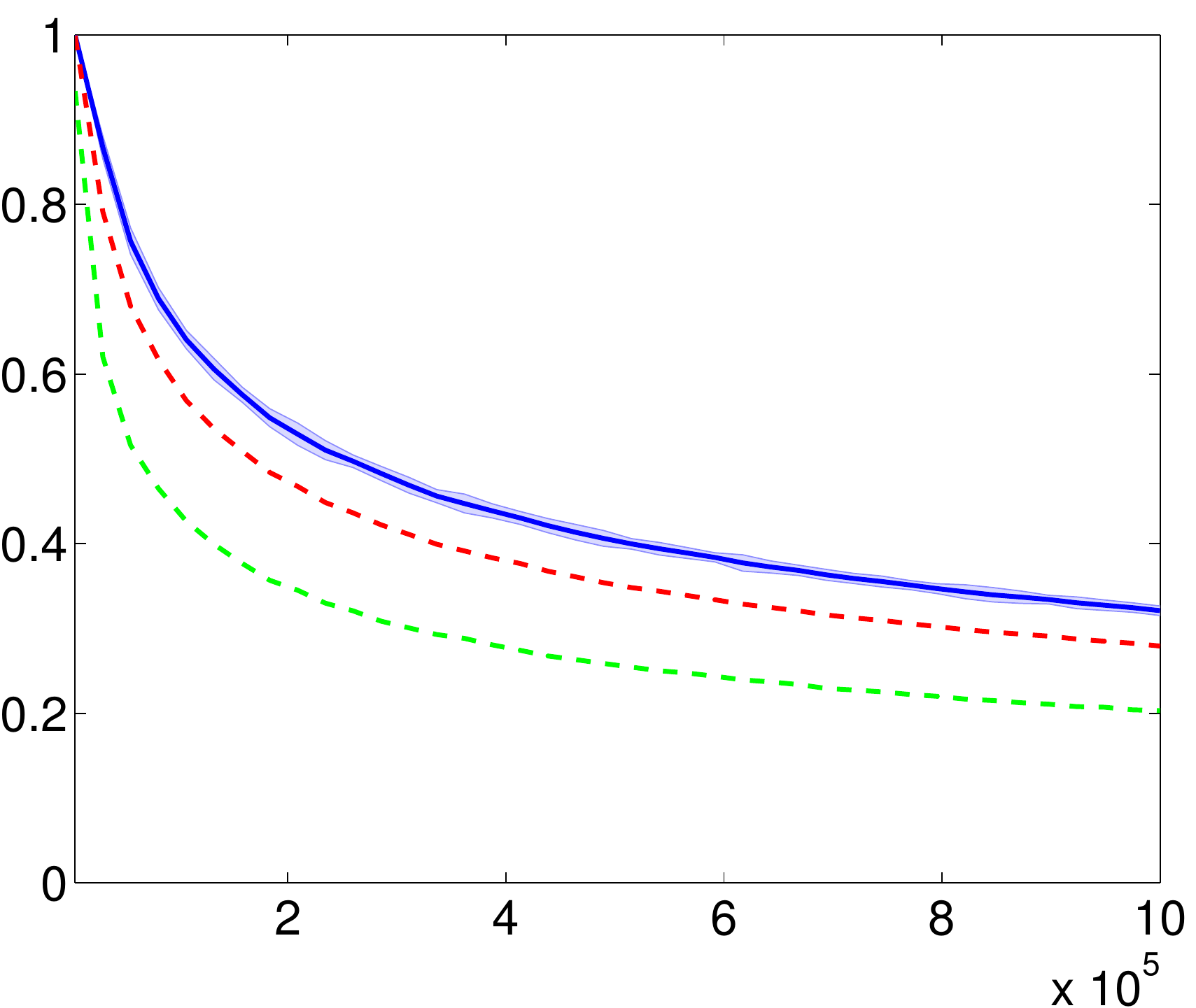}}
        \subfigure[Schr\"{o}dinger cat $q_0=3$]{\includegraphics[scale=0.42]{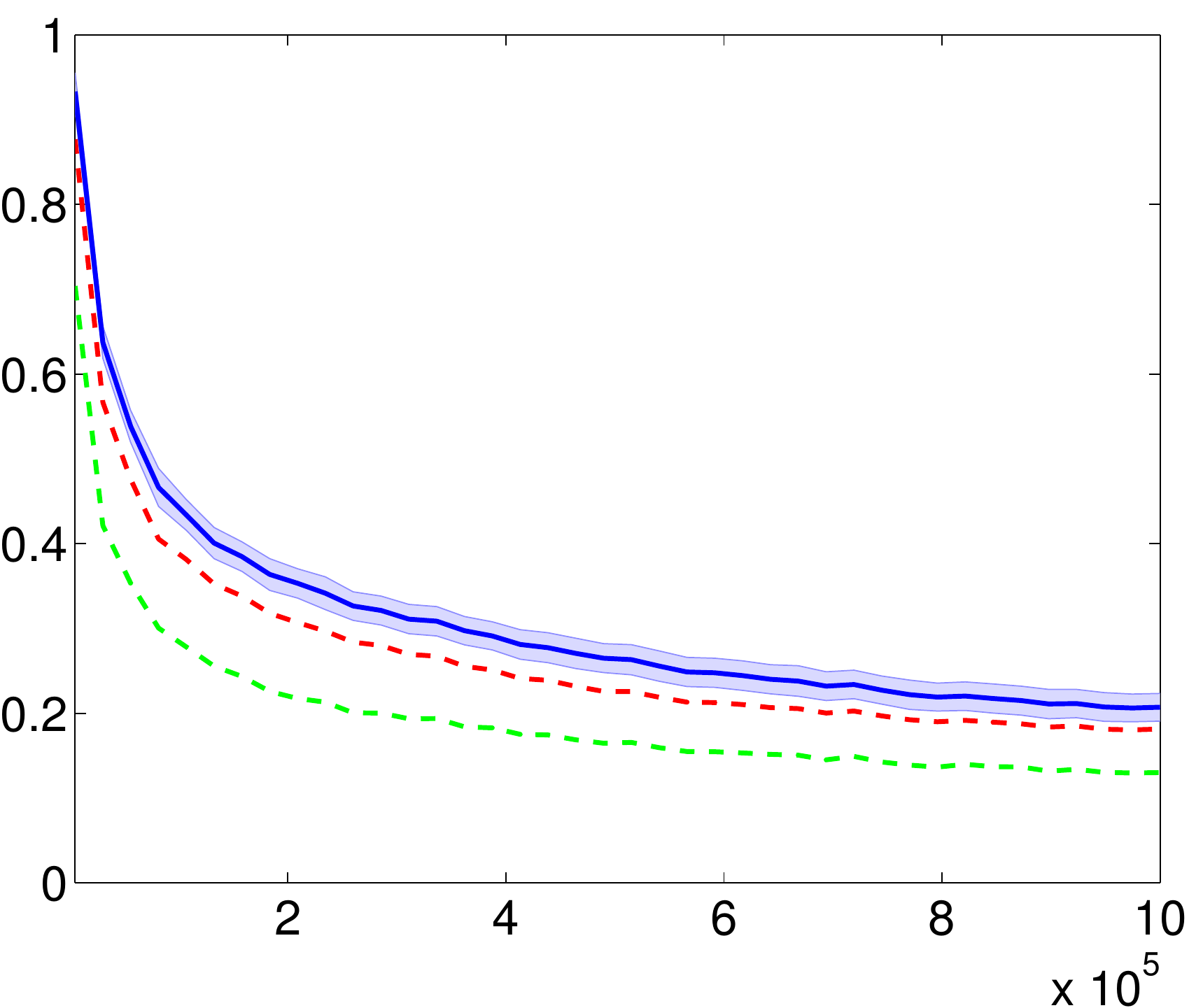} }
        \subfigure[Thermal $\beta=1/10$]{\includegraphics[scale=0.42]{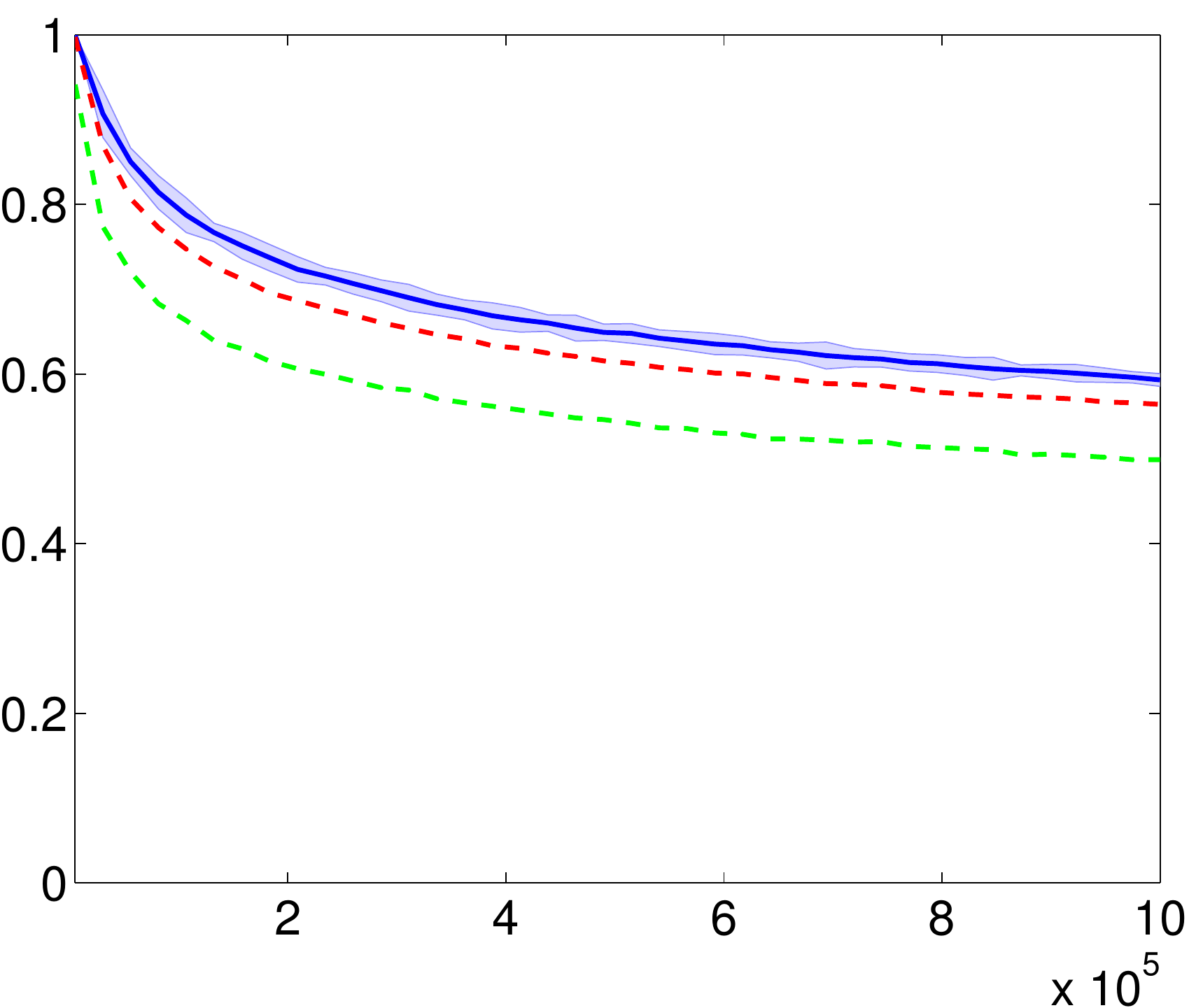} }
        \subfigure[Thermal $\beta=1/4$]{\includegraphics[scale=0.42]{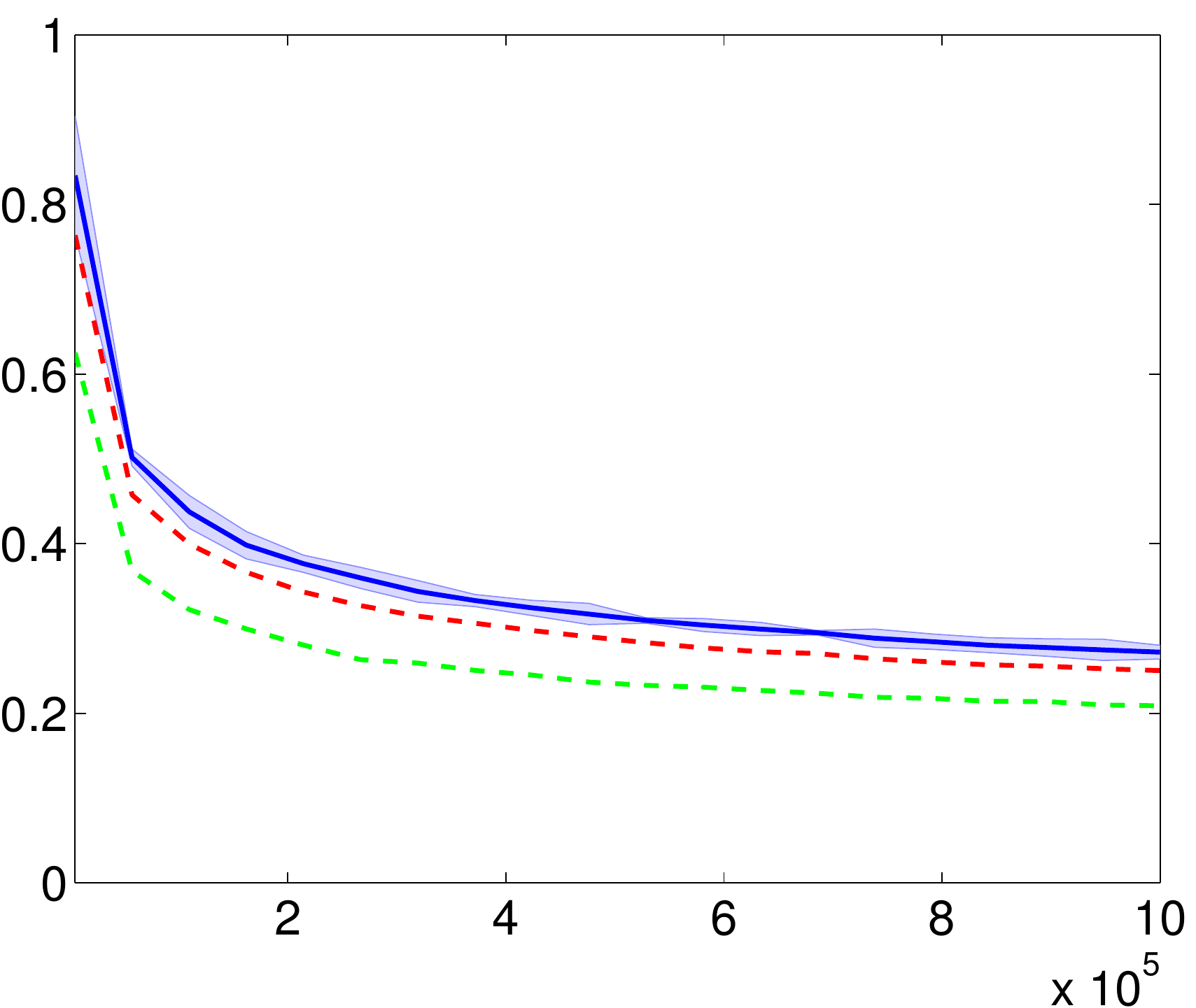} }
\end{figure}

We estimate numerically the (relative) root mean square error (RMSE) 
$$
	\textrm{RMSE}(n) =  \| \tilde \rho^\eta - \rho \|_2 / \| \rho \|_2 
$$
of our soft thresholding estimator. More precisely, Figure~\ref{fig:RMSE} shows the evolution with $n$ of the expected value of the RMSE. This expected value is evaluated by an empirical mean with Monte Carlo simulation using 50 replications for each value of $n$. To evaluate the deviation with respect to this mean, we also display the confidence interval at $\pm 3$ times the standard deviation of the RMSE. 

\noindent The threshold values $t_{j,k}$ that are used in (\ref{rho-seuil}) to define our estimator are somewhat conservative. In practice, smaller values offer better decay of the RMSE. Figure \ref{fig:RMSE} displays in dashed red (resp. dashed green) the decay of the RMSE obtained using thresholds $0.8 t_{j,k}$ (resp. $0.5 t_{j,k}$). We found on these three examples and for $\eta=0.9$ that using $0.5 t_{j,k}$ gives consistently the lowest RMSE among other choices of thresholds proportional to the $t_{j,k}$ values.\\

\noindent We found numerically that the decay of the RMSE with $n$ almost perfectly fits a power-law, which (up to logarithmic factor) is in accordance with the upper-bounds of Corollary~\ref{coro1}. Following this Corollary in the setting $\eta\in(\frac{1}{2},1)$ and $ r=2$,  we fit a power law of the form 
$$
	\mathbb{E}(\textrm{RMSE}(n)) \approx n^{ -\frac{\tilde B}{2 ( 4 \gamma+\tilde B)} }.
$$
We perform a linear regression in a log-log domain to estimate $\tilde B$. 
Table~\ref{tab:B} reports the estimated value of $\tilde B$ we found using this procedure. 

\begin{table}[!h]
\caption{Estimated values of $\tilde B$ when using $\eta=0.9$, $\varepsilon=1$ and $N=30$.}
\label{tab:B}
\begin{center}
\begin{tabular}{|c|c|c|c|}
\hline 
Coherent $q_0=3$ &
Schr\"{o}dinger cat $q_0=3$ &
Thermal $\beta=1/10$ & 
Thermal $\beta=1/4$ \\ \hline
$\tilde B \approx 0.174$ &
$\tilde B \approx 0.227$ &
$\tilde B \approx 0.037$ & 
$\tilde B \approx 0.082$ \\\hline
\end{tabular}
\end{center}
\end{table}

\bigskip
\bigskip
\bigskip
\ack  We would like to thank Arnak Dalalyan for his relevant remarks and for his rereading of the present paper.  The research of Katia Meziani is partly supported by the french Agence Nationale de la Recherche (ANR 2011 BS01 010 01 projet Calibration )


\newpage
\appendix

\section{Proof of Proposition \ref{thm.oracle}}
\label{sec.proofs.th}

\noindent The proofs follow the main lines of \cite{Alquier2008a,Alquier2008b}. First, we
need
a set of preliminary lemmas.

\subsection{Some preliminary results}
First, we remind Hoeffdig's inequality for bounded random variables.
\begin{lemma}
Let us assume that $Z_1$, ..., $Z_n$ are independent real-valued random variables with
$a_i\leq |Z_i|\leq b_i$. Then, for any $\lambda>0$,
$$ \mathbb{P}\left(\left|\sum_{i=1}^{n} \left[Z_i - \mathbb{E}(Z_i)\right]\right|
        \geq \lambda \right) \leq 2 \exp\left( -\frac{2 \lambda^2}{\sum_{i=1}^{n}(b_i-a_i)^2}
        \right).$$
\end{lemma}

As a consequence, we have the following inequality for complex random variables.
\begin{lemma} \label{lemme_hoeffding_complexe}
Let us assume that $Z_1$, ..., $Z_n$ are independent complex-valued random variables with
$|Z_i|\leq c$. Then, for any $t>0$,
$$
  \mathbb{P}\left(\left|\frac{1}{n}\sum_{i=1}^{n} \left[Z_i - \mathbb{E}(Z_i)\right]\right|
        \geq t \right) 
        \leq 4 \exp\left(-\frac{n t^2}{4c^2}\right).
        $$
\end{lemma}

\noindent \textit{Proof: } We have
\begin{eqnarray*}
  \mathbb{P}\left(\left|\sum_{i=1}^{n} \left[Z_i - \mathbb{E}(Z_i)\right]\right|
        \geq \lambda \right) & \leq
  \mathbb{P}\left(\left|{\rm Re}\sum_{i=1}^{n} \left[Z_i - \mathbb{E}(Z_i)\right]\right|
        \geq \frac{ \lambda}{\sqrt{2}} \right)
       \\
       & \hspace{1cm} +
   \mathbb{P}\left(\left| {\rm Im}\sum_{i=1}^{n}\left[Z_i - \mathbb{E}(Z_i)\right]\right|
        \geq  \frac{ \lambda}{\sqrt{2}}  \right) \\
 & \leq  \mathbb{P}\left(\left|\sum_{i=1}^{n} \left[{\rm Re}(Z_i) - \mathbb{E}({\rm Re}(Z_i))
             \right]\right|
        \geq \frac{ \lambda}{\sqrt{2}} \right)
       \\
       & \hspace{1cm} +
   \mathbb{P}\left(\left|\sum_{i=1}^{n}\left[ {\rm Im}(Z_i) - \mathbb{E}( {\rm Im}(Z_i))
    \right]\right|
        \geq  \frac{ \lambda}{\sqrt{2}}  \right).
\end{eqnarray*}
Now, we apply Hoeffding's inequality to the random variables ${\rm Re}(Z_i)$ which satisfy
$-c\leq {\rm Re}(Z_i) \leq c$. So we have:
$$ \mathbb{P}\left(\left|\sum_{i=1}^{n} \left[{\rm Re}(Z_i) - \mathbb{E}({\rm Re}(Z_i))
             \right]\right|
        \geq \frac{ \lambda}{\sqrt{2}} \right)
         \leq 2 \exp\left(-\frac{2 \left(\frac{\lambda}{\sqrt{2}}\right)^2}
              {\sum_{i=1}^{n}(2c)^2}\right)
               = 2 \exp\left(-\frac{\lambda^2}{4nc^2}\right). $$
We have exactly the same result for ${\rm Im}(Z_i)$ so finally:
$$
  \mathbb{P}\left(\left|\sum_{i=1}^{n} \left[Z_i - \mathbb{E}(Z_i)\right]\right|
        \geq \lambda \right) 
        \leq 4 \exp\left(-\frac{\lambda^2}{4nc^2}\right).
        $$
Put $t= \lambda/n$ to end the proof. $ \Box $

\begin{lemma}
                \label{lem.deviation}
                For some fixed $\varepsilon\in(0,1)$, let us define the set
                $$
                \Omega_{\varepsilon}\deq \left\{ \forall (j,k)\in J(N),\quad
\left|\hat{\rho}^\eta_{j,k} - \rho_{j,k} \right| \leq t_{j,k} \right\},
                $$
                where the $( t_{j,k})_{j,k}$ are defined in (\ref{tjk}) and the
set $J(N)$ is defined in (\ref{eqJ}). Then we have
                $$P(\Omega_{\varepsilon})\geq1-\varepsilon.$$
\end{lemma}

\noindent \textit{Proof: }Lemma~\ref{lem.deviation} is proved by using  Hoeffding's inequality. In
this aim, we have to first notice
$$E_\rho[\hat{\rho}^\eta_{j,k}]=\rho_{j,k}.$$
Indeed, by using (\ref{Gjk}) , (\ref{fourierproun}), (\ref{eq:patterneta}) and
(\ref{rhojk}), we have 
\begin{eqnarray*}
        E_\rho[\hat{\rho}^\eta_{j,k}] &=&
        E_\rho[G_{j,k}(\frac{Y}{\sqrt\eta},\Phi)] =
        E_\rho[f_{j,k}^{\eta}(\frac{Y}{\sqrt\eta})e^{-i(j-k)\Phi}]
         \\
        & =&\frac{1}{\pi}\int_0^\pi e^{-i(j-k)\phi}\int
        f_{j,k}^{\eta}(y)\sqrt{\eta}p_\rho^\eta(y\sqrt{\eta}|\phi)dy
        d\phi\\
        &= &\frac{1}{\pi}\int_0^\pi e^{-i(j-k)\phi}\frac{1}{2\pi}\int
        \widetilde{f}_{j,k}^{\eta}(t)
        \mathcal{F}_1[\sqrt{\eta}p_\rho^\eta(\cdot\sqrt{\eta}|\phi)](t)dt
        d\phi\\
        &=& \frac{1}{\pi} \int_0^\pi e^{-i(j-k)\phi} \frac{1}{2\pi}
        \int \widetilde{f}_{j,k}(t) e^{\gamma
          t^2} \mathcal{F}_1[p_\rho(\cdot|\phi)](t)
        \widetilde{N}^\eta(t) dt d\phi\\
                &=& \frac{1}{\pi}    \int_0^\pi   \int  e^{-i(j-k)\phi} 
    f_{j,k}(x) p_\rho(x|\phi)dx d\phi=\rho_{j,k}.
\end{eqnarray*} 
 Moreover, we easily get from the definition of $G_{j,k}$ in (\ref{Gjk}) that
for all $\ell=1\cdots,n$ and  $\forall (j,k)\in J(N)$
   $$ 
\left| G_{j,k}\paren{\frac{Y_\ell}{\sqrt\eta},\Phi_\ell}
\right| \leq\|f_{j,k}^{\eta}\|_{\infty}.
   $$
Then, for 
   $$
   t_{j,k} =
2 \|f_{j,k}^{\eta}\|_{\infty}\sqrt{\frac{\log\left(\frac{2 N(N+1)}{\varepsilon}
\right)}{n}}
  $$
  and according to Lemma~\ref{lemme_hoeffding_complexe}
$$
P \left( \left|   \hat{\rho}^\eta_{j,k} - \rho_{j,k} \right| \geq t_{j,k}   \right)
\leq 4 \exp\left[ - \frac{nt_{j,k}^{2}}{4
\|f_{j,k}^{\eta}\|_{\infty}^{2}}\right]
          = \frac{2 \varepsilon}{N(N+1)} . 
$$
By the classical union bound argument:
 \begin{equation*}
          P \left(\Omega_{\varepsilon}^c  \right) 
          \leq \sum_{ (j,k)\in J(N)} P \left( \left|  \hat{\rho}^\eta_{j,k} -
\rho_{j,k}
         \right| \geq t_{j,k}    \right)\leq \sum_{ (j,k)\in J(N)} \frac{2
\varepsilon}{N(N+1)}\leq \varepsilon .
\end{equation*}$ \Box $

\begin{lemma}
           \label{lem.projection}
            For some fixed $\varepsilon\in(0,1)$ and $\forall(j,k)\in J(N)$,
with $J(N)$ defined in (\ref{eqJ}), we define the set
            $$
             R_{j,k}^{\varepsilon}\deq \left\{ \nu \,\,\mathrm{ density}\,\,\mathrm{  matrix},
|\nu_{j,k} - \hat{\rho}^\eta_{j,k}| \leq t_{j,k} \right\},
            $$
            where the $( t_{j,k})_{j,k}$ are defined in (\ref{tjk}). Then, on
the event $\Omega_{\varepsilon}$ defined in Lemma~\ref{lem.deviation} and
$\forall (j,k)\in J(N)$
          \begin{enumerate}
             \item $\rho\in R_{j,k}^{\varepsilon}$.
             \item $R_{j,k}^{\varepsilon}$ is closed and convex set.
             \item For $\Pi_{j,k}^{\varepsilon} $ the orthogonal projection onto
$ R_{j,k}^{\varepsilon}$ and for any density matrix $\nu$,
                 \begin{equation}
                   \label{eq.improv}
                   \left\|\rho - \Pi_{j,k}^{\varepsilon}(\nu) \right\|^{2}_{2} 
\leq \left\|\rho - \nu \right\|^{2}_{2} .
                 \end{equation}
         \end{enumerate}
\end{lemma}

\noindent \textit{Proof: }The first point is just a consequence of Lemma~\ref{lem.deviation}. The second point  comes from the definition of $R_{j,k}^{\varepsilon}$. \\
Moreover, it is well known that for any closed and convex set $\mathcal{C}$, if
$\Pi_{\mathcal{C}}$ is the orthogonal projection on $\mathcal{C}$, the following
property holds:
$$ 
\forall x\in\mathcal{C} ,\forall y,\quad \|\Pi_{\mathcal{C}}(y) - x \|_2 \leq
\|y-x\|_2 . 
$$
This concludes the proof of the third point. $ \Box $

\begin{lemma}
          \label{lem.formule}
          For $\varepsilon\in(0,1)$, any fixed $(j,k)\in J(N)$, with $J(N)$
defined in (\ref{eqJ}), and any density matrix $\nu$, we denote by $\nu'$ the
projection of $\nu$ into $R_{j,k}^{\varepsilon}$,
          $$
          \nu' \deq \Pi_{j,k}^{\varepsilon}(\nu)=[\nu'_{\ell,m}]_{\ell,m},
          $$
          with $R_{j,k}^{\varepsilon}$ defined in Lemma~\ref{lem.projection}.
Then, the entries $\nu'_{\ell,m}$  of $\nu'$ are equal to
          $$
          \nu'_{\ell,m} =\left\{
          \begin{array}{l}
          \nu_{j,k} + \frac{\hat{\rho}_{j,k}^{\eta} - \nu_{j,k}}
                   {\left|\hat{\rho}_{j,k}^{\eta} - \nu_{j,k}\right|}
\left(\left|\hat{\rho}_{j,k}^{\eta} -
\nu_{j,k}\right|-t_{j,k}\right)_{+},\mathrm{ if } (\ell,m)=(j,k),\\\\
          \nu_{\ell,m}, \mathrm{ otherwise,}
           \end{array}\right.
           $$
           with the convention $0/0=0$.
\end{lemma}

\begin{center}
\includegraphics[width=8cm]{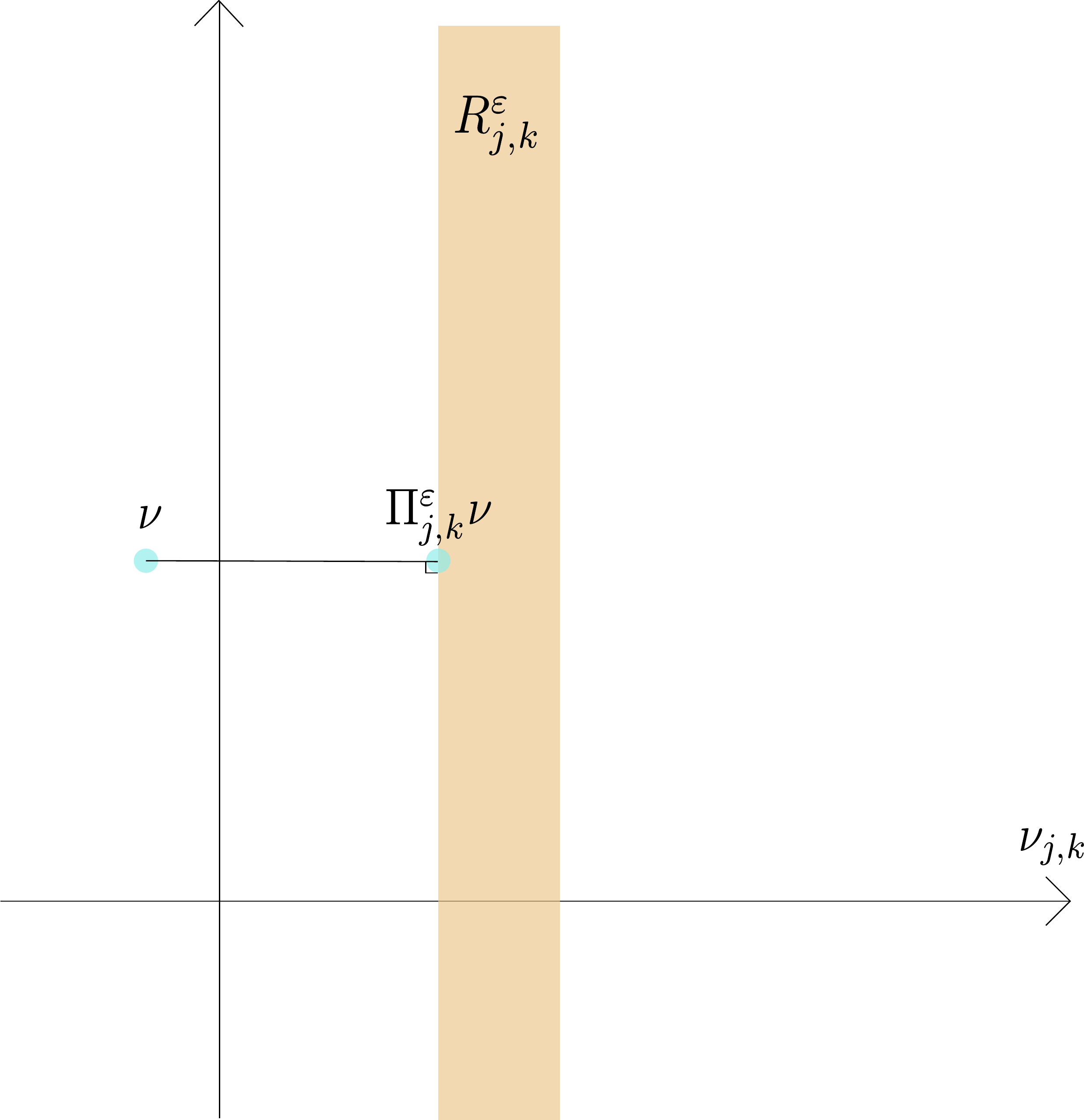}
\end{center}

\noindent \textit{Proof: }The projection $\nu'$ of $\nu$ into $R_{j,k}^{\varepsilon}$ satisfies
      $$
      \nu' = \arg\min_{x\in R^{\varepsilon}_{j,k}} \|\nu - x\|_{2}^{2}  =
\arg\min_{x\in R^{\varepsilon}_{j,k}} \sum_{\ell,m=0}^{\infty} \left|x_{\ell,m}
- \nu_{\ell,m} \right|^{2} .
      $$
     As the constraint $x\in R^{\varepsilon}_{j,k}$ is only a constraint on
$x_{j,k}$, it is clear that for $(\ell,m)\neq(j,k)$ the minimum is reached for
$x_{j,k}=\nu_{j,k}$. Finally,
    $$
    \nu'_{j,k} = \arg\min_{x_{j,k}:\,|x_{j,k} - \hat{\rho}_{j,k}^{\eta}|\leq
t_{j,k} } |\nu_{j,k} - x_{j,k}|^{2}. 
    $$
The solution $\nu'_{j,k}$ is obvious:
             \begin{eqnarray*}
                \nu'_{j,k} & = &\left\{\begin{array}{l l}
                  \nu_{j,k} & \mathrm{ if }\, \,\,|\nu_{j,k} - \hat{\rho}_{j,k}^{\eta}|\leq
t_{j,k},\\\\
                \nu_{j,k} + \frac{\hat{\rho}_{j,k}^{\eta} - \nu_{j,k}}
                   {\left|\hat{\rho}_{j,k}^{\eta} - \nu_{j,k}\right|}
\left(\left|\hat{\rho}_{j,k}^{\eta} -
\nu_{j,k}\right|-t_{j,k}\right) & \mathrm{ otherwise}
                   \end{array}\right.\\
                & = &\nu_{j,k} + \frac{\hat{\rho}_{j,k}^{\eta} - \nu_{j,k}}
                   {\left|\hat{\rho}_{j,k}^{\eta} - \nu_{j,k}\right|}
\left(\left|\hat{\rho}_{j,k}^{\eta} -
\nu_{j,k}\right|-t_{j,k}\right)_{+}.
            \end{eqnarray*}
This ends the proof. $ \Box $

\begin{dfn} 
    For $m>0$ an integer, let 
    $$
    I \deq \{(j_1,k_1),\dots,(j_m,k_m)\}\subseteq J(N)
    $$
    be a set of indices, where $J(N)$ is the set defined in (\ref{eqJ}). 
that $\forall \ell\neq i$, $(j_\ell,k_\ell)\neq (j_i,k_i)$. 
      For $\varepsilon\in(0,1)$ and for any density matrix $\nu$, we denote by
$\Pi_{I}^{\varepsilon}(\nu)$ the successive projections of $\nu$ into spaces
$\left(R_{j_i,k_i}^{\varepsilon}\right)_{j_i,k_i}$, i.e.
      $$
       \Pi_{I}^{\varepsilon}(\nu)\deq\Pi_{j_m,k_m}^{\varepsilon}
\Pi_{j_{m-1},k_{m-1}}^{\varepsilon} \ldots \Pi_{j_2,k_2}^{\varepsilon}
\Pi_{j_1,k_1}^{\varepsilon}( \nu ).
     $$
  \end{dfn}
Note that  for any set of indices $I$ and from Lemma~\ref{lem.formule}, the
application of the successive projections $\Pi_{I}^{\varepsilon}$ to a density
matrix $\nu$
does not depend on the order of the successive projections.

\begin{lemma}
         \label{lem.estimator}
         For $\varepsilon\in(0,1)$,  for $J(N)$ defined in (\ref{eqJ}) and for
$\tilde{\rho}^{\eta}$ defined in (\ref{rho-seuil}), we have
          $$
          \tilde{\rho}^\eta = \Pi_{J(N)}^{\varepsilon}({\bf 0}),
          $$
          where ${\bf 0}$ is the zero-infinite matrix.
\end{lemma}

\noindent \textit{Proof: }This is obvious from the definition of $\tilde{\rho}^\eta$ and from
Lemma~\ref{lem.formule} applied to $\nu={\bf 0}$. $ \Box $

\subsection{Proof of Proposition \ref{thm.oracle}}

\noindent \textit{Proof: }For $J(N)$ the set of indices defined in (\ref{eqJ}), let $I$ be a subset of
$J(N)$, $I\subseteq J(N)$. For a fixed $\varepsilon\in(0,1)$, we have by
Lemma~\ref{lem.estimator} and by successive applications of  the
inequality(\ref{eq.improv}) to all pair of indices $(j,k)\notin I$

 \begin{equation}
 \label{eq1}
             \| \tilde{\rho}^\eta-\rho \|_{2}^{2} = \| \Pi_{J }^{\varepsilon}(
{ \bf 0} )-\rho \|_{2}^{2}  \leq \| \Pi_{I }^{\varepsilon}( {\bf 0} )-\rho
\|_{2}^{2}.
  \end{equation}
Moreover, from Lemma~\ref{lem.formule} applied to $\nu={\bf 0}$, we get 
$$
(\Pi_{I}^{\varepsilon} ({\bf 0}))_{j,k}= \left\{
\begin{array}{l l} \frac{\hat{\rho}_{j,k}^{\eta} }
                   {\left|\hat{\rho}_{j,k}^{\eta}\right|}
\left|\left|\hat{\rho}_{j,k}^{\eta}\right|-t_{j,k}\right)_{+} & \mathrm{ if } (j,k)\in I,\\
0 & \mathrm{ otherwise.}\end{array}\right.
$$
Therefore, from (\ref{eq1}) we get
             \begin{eqnarray}
             \label{eq2}
               \| \tilde{\rho}^\eta-\rho \|_{2}^{2}
             & \leq& \sum_{j,k=0}^{\infty} |\rho_{j,k} -
(\Pi_{I}^{\varepsilon} ({\bf 0}))_{j,k} |^{2}\nonumber\\
             & =& \sum_{(j,k)\in I} \left|\rho_{j,k} -  \frac{\hat{\rho}_{j,k}^{\eta} }
                   {\left|\hat{\rho}_{j,k}^{\eta}\right|}
\left(\left|\hat{\rho}_{j,k}^{\eta}\right|-t_{j,k}\right)_{+}
\right|^{2} + \sum_{(j,k)\notin I} |\rho_{j,k}|^{2}\nonumber\\
             &\deq &\sum_{(j,k)\in I} |A_{j,k}|^{2} +
\sum_{(j,k)\notin I} |\rho_{j,k}|^{2},
             \end{eqnarray}
where
\begin{eqnarray*}
A_{j,k}&=&\rho_{j,k} - \frac{\hat{\rho}_{j,k}^{\eta} }
                   {\left|\hat{\rho}_{j,k}^{\eta}\right|}
\left(\left|\hat{\rho}_{j,k}^{\eta}\right|-t_{j,k}\right)_{+}\\
&=& \left\{ \begin{array}{ccc}
                       \rho_{j,k}, & \mathrm{if  }|\hat{\rho}_{j,k}^{\eta}|\leq
t_{j,k},\\\\
                          \rho_{j,k}-\frac{\hat{\rho}_{j,k}^{\eta} }
                   {\left|\hat{\rho}_{j,k}^{\eta}\right|}
\left(\left|\hat{\rho}_{j,k}^{\eta}\right|-t_{j,k}\right),& \mathrm{ otherwise}.
                             \end{array} \right.
\end{eqnarray*}
Moreover
\begin{eqnarray*}
|A_{j,k}|& \leq & \left\{ \begin{array}{ccc}
                     \left|\hat{\rho}^\eta_{j,k} - \rho_{j,k} \right|
+|\hat{\rho}_{j,k}^{\eta}|, & \mathrm{if  }|\hat{\rho}_{j,k}^{\eta}|\leq
t_{j,k},\\\\
                      \left|\hat{\rho}^\eta_{j,k} - \rho_{j,k} \right|+
t_{j,k},& \mathrm{otherwise}
                             \end{array} \right.\\
                             &\leq& \left|\hat{\rho}^\eta_{j,k} - \rho_{j,k}
\right|+ t_{j,k}.
\end{eqnarray*}

For any $(j,k)\in I$ and on the event $\Omega_{\varepsilon}$  defined in
Lemma~\ref{lem.deviation}, it holds
$$
\left|\hat{\rho}^\eta_{j,k} - \rho_{j,k} \right| \leq t_{j,k}.
$$
Therefore from (\ref{eq2}) 
              \begin{equation*}
               \| \tilde{\rho}^\eta-\rho \|_{2}^{2} \leq \sum_{(j,k)\in I} ( 2
t_{j,k} )^{2} + \sum_{(j,k)\notin I}| \rho_{j,k}|^{2}.
             \end{equation*}
We conclude the proof by taking the infimum over the set $I\subseteq J(N)$. $ \Box $

\section{Proof of Theorem \ref{coro1}}
\label{sec.proofs.co1}

\noindent \textit{Proof: }For $r_{0}\in(0,2)$, $B_{0}>0$ and $N$ as in (\ref{Ncorr1}), let $M$ be an integer s.t. $M<N$. We define the set 
$$J(M) := \{(j,k)\in \mathbb{N}^2,\, 0\leq j+k\leq M\}.$$
 Then, for $\varepsilon\in (0,1)$ and by applying Proposition~\ref{thm.oracle} to $I
= J(M)$, with probability larger than $1-\varepsilon$, we obtain
       \begin{eqnarray}
               \nonumber
                \left\|\tilde{\rho}^\eta - \rho \right\|^{2}_{2} 
                & \leq& \inf_{0\leq M \leq N-1}\left\{ 4 \sum_{0\leq j+k\leq M}
t_{j,k}^{2} + \sum_{j+k> M}|\rho_{j,k}|^{2}  \right\} \\
                & = &\inf_{0\leq M \leq N-1}\left\{ \frac{16}{n} \sum_{0\leq
j+k\leq M} \|f_{j,k}^\eta\|_{\infty}^{2}\log\left(2N(N+1)/\varepsilon\right)+\right.\nonumber\\
&&\left.
\sum_{j+k> M}|\rho_{j,k}|^{2} \right\}.
             \label{etape1}
     \end{eqnarray}
 \noindent\\
     
\paragraph{\textit{a) For $\eta=1$ and $r\in[r_0,2]$.}}
\noindent\\

 %
As  $ f_{j,k}^\eta=f_{j,k}$ for $\eta=1$, we have by pluging  (\ref{Katia1}) and
(\ref{normf}) into (\ref{etape1}) 
       \begin{eqnarray}
               \nonumber
                \left\|\tilde{\rho}^\eta - \rho \right\|^{2}_{2} 
                & \leq &\inf_{0\leq M \leq N-1}\left\{ \frac{16}{n} \sum_{0\leq
j+k\leq M} \|f_{j,k}\|_{\infty}^{2}\log\left(2N(N+1)/\varepsilon\right)+
\sum_{j+k> M}|\rho_{j,k}|^{2} \right\}.
                 \nonumber\\
                 &\leq& \inf_{0\leq M \leq N-1}\left\{\frac{c_1}{n}
M^{\frac{10}{3}}\log(N/\varepsilon)  + \mathcal{C} M^{2-\frac{r}{2}}
e^{-2BM^{\frac{r}{2}}} \right\},
 \label{etape2}
     \end{eqnarray}
for some constant $c_1>0$.\\
 For $N$ such in (\ref{Ncorr1}) and by taking $M=(\log(n)/2B)^{2/r}<N$, it leads
to
$$
 \left\|\tilde{\rho}^\eta  - \rho \right\|^{2}_{2}  \leq\mathcal{C}_1 \log
\left(\log(n)/\varepsilon\right)
 \left(\log(n)\right)^{\frac{20}{3r}} n^{-1}
$$
for some constant $\mathcal{C}_{1}>0$.
\noindent\\

\paragraph{\textit{b) For $\eta\in(1/2,1)$ and $r=2$.}}
\noindent\\

%
Next, we deal with the case $1>\eta>1/2$. We plug (\ref{Katia1}) and
(\ref{normfeta}) into (\ref{etape1}) to obtain in the case $r=2$
              \begin{equation*}
                        \left\|\tilde{\rho}^\eta  - \rho \right\|^{2}_{2}
\leq\inf_{0\leq M <N }\left\{ \frac{c_{2}}{n}  \log\left(N/
                        \varepsilon\right)M^{\frac{1}{3}} e^{8 \gamma M}+
\mathcal{C}M  e^{-2BM}\right\},
               \end{equation*}
for some constant $c_2>0$.\\
By taking $M=M(n)$ s.t.
$$ 
M=\frac{\log(n)}{2(4\gamma+B)} ,
$$
we obtain
$$
 \left\|\tilde{\rho}^\eta - \rho \right\|^{2}_{2}  \leq\mathcal{C}_2  
n^{-\frac{B}{4\gamma+B}}\left(\log\left( \log (n)/\varepsilon\right)
\left(\log(n)\right)^{1/3}+\log(n) \right) ,
$$
for some constant $\mathcal{C}_2 >0$. 
\noindent\\ 

\paragraph{\textit{c) For $\eta\in(1/2,1)$ and $r\in(r_0,2)$.}}
\noindent\\
Finally, in the case $\eta\in(1/2,1)$ and $r\in(r_0,2)$ and by plugging 
(\ref{Katia1}) and (\ref{normfeta}) into (\ref{etape1}) we get:
              \begin{equation*}
                        \left\|\tilde{\rho}^\eta  - \rho \right\|^{2}_{2}
\leq\inf_{0\leq M <N }\left\{ \frac{c_3}{n}\log\left(N
                        /\varepsilon\right) M^{\frac{1}{3}} e^{8 \gamma M}+
\mathcal{C} M^{2-\frac{r}{2}}  e^{-2BM^{\frac{r}{2}}}\right\},
               \end{equation*}
for some constant $c_3>0$.\\
For $M$ a solution of the equation $8\gamma M +2 B M^{\frac{r}{2}} = \log(n)$
and for 
$$ M(n)= \frac{1}{8\gamma} \log(n) - \frac{2B}{(8\gamma)^{1+r/2}} \log(n)^{r/2}
+ o(\log(n)^{r/2})$$
in particular, we obtain 
$$
 \left\|\tilde{\rho}^\eta - \rho \right\|^{2}_{2}  \leq\mathcal{C}_3  
\exp^{-2BM^{r/2}} \left( \log(n)^{2-r/2} +
 \log(n)^{1/3}\log\left(N/\varepsilon\right) \right),
$$
for some constant $\mathcal{C}_3 >0$. $ \Box $

\section{Proof of Theorem \ref{coro2}}
\label{sec.proofs.co2}

\noindent\textit{Proof: }We apply Theorem \ref{thm.oracle} for $I = \{(j_0,j_0)\}$. We obtain, with
probability larger than
$1-\varepsilon$,
           \begin{eqnarray*}
                      \left\|\tilde{\rho}^\eta - \rho \right\|^{2}_{2}& \leq&
\frac{16}{n} \sum_{(j,k)=(j_0,j_0)} \|f_{j,k}\|_{\infty}^{2}
                     \log\left(2N(N+1)/\varepsilon\right)+ \sum_{(j,k)\neq
(j_0,j_0)}\rho_{j,k}^{2}  \\
                       & =& \frac{16}{n} \|f_{j_0,j_0}\|_{\infty}^{2}
                     \log\left(2N(N+1)/\varepsilon\right)+ 0.
           \end{eqnarray*}
For $n$ large enough,  $N=N(n)\geq2$. Then, $(N+1)< 2N$ and 
           \begin{eqnarray*}
                      \left\|\tilde{\rho}^\eta - \rho \right\|^{2}_{2}& < &
\frac{16}{n} \|f_{j_0,j_0}\|_{\infty}^{2}
                     \log\left(4N^2/\varepsilon\right) \\
                     & =& \frac{16}{n} \|f_{j_0,j_0}\|_{\infty}^{2} \left[
                     2 \log(N) + \log\left(4/\varepsilon\right) \right] \\
                     & \leq & \frac{16}{n} \|f_{j_0,j_0}\|_{\infty}^{2} \left[
                     \frac{4}{r_{0}} \log\left(\frac{\log(n)}{2B_{0}}\right) +
\log\left(4/\varepsilon\right) \right]
           \end{eqnarray*}
where we replaced $N$ by its definition. As $r_0 < 2$, $4/r_0 > 1$ and we have
the following rough upper bound:
 \begin{eqnarray*}
                      \left\|\tilde{\rho}^\eta - \rho \right\|^{2}_{2} & <&
\frac{16}{n} \|f_{j_0,j_0}\|_{\infty}^{2} \left[
                     \frac{4}{r_{0}} \log\left(\frac{\log(n)}{2B_{0}}\right) + 
\frac{4}{r_{0}}  \log\left(4/\varepsilon\right) \right] \\
                       & =& \frac{64}{nr_{0}} \|f_{j_0,j_0}\|_{\infty}^{2}
\log\left(\frac{2\log(n)}{B_{0} \varepsilon}\right).
           \end{eqnarray*} $ \Box $

\section{Technical Lemmas}
\label{sec.annexe}
\noindent Useful lemmas are the following. 

\begin{lemma}
         \label{biais}
         For $\rho\in\mathcal{R}(C,B,r)$, the set defined in
(\ref{eq.classcoeff}), there exists a $M_0$ s.t. $\forall M\geq M_0$   implies
                    \begin{equation}
                     \label{Katia1}
                     \sum_{j+k> M}|\rho_{j,k}|^{2} \leq \mathcal{C}
M^{2-\frac{r}{2}} e^{-2BM^{\frac{r}{2}}},
                  \end{equation}
                     where $\mathcal{C}=\frac{2C^2}{Br}.$
\end{lemma}

\noindent \textit{Proof: }For $\rho\in\mathcal{R}(C,B,r)$, we have by the definition of the class
$\mathcal{R}(C,B,r)$ and by Lemma~3 in \cite{ABM}
                  \begin{eqnarray*}
                  \sum_{j+k> M}|\rho_{j,k}|^{2}\leq C^2  \sum_{j+k> M}\exp(-2B
(j+k)^{r/2})\leq  \frac{2C^2}{Br} M^{2-\frac{r}{2}} e^{-2BM^{\frac{r}{2}}}.
                 \end{eqnarray*} $ \Box $

\begin{lemma}
         \label{norm}
         For $\eta\in(1/2,1)$, there exists a positive constant
$C^\eta_{\infty}>0$ s.t.
          \begin{equation}
                         \label{normfeta}
                         \sum_{0\leq j+k\leq M} \|f_{j,k}^{\eta}\|_{\infty}^{2}
\leq C_{\infty}^{\eta} M^{\frac{1}{3}}  e^{8\gamma M},
          \end{equation}
          where $\gamma=(1-\eta)/(4\eta)$ and the $(f_{j,k}^{\eta})_{j,k}$ are
the adapted pattern functions defined in expression (\ref{eq:patterneta}).\\
          There exists a positive constant $C_{\infty}>0$ s.t.
           \begin{equation}
                   \label{normf}
                   \sum_{0\leq j+k\leq M} \|f_{j,k}\|_{\infty}^{2} \leq C_\infty
M^{\frac{10}{3}},
           \end{equation}
          where the $ (f_{j,k})_{j,k}$ are the pattern functions defined in
expression (\ref{Rpattern}).
 \end{lemma}

\noindent \textit{Proof: }For the proof of this lemma, we refer to Lemma~4 and Lemma~5 in \cite{ABM}.  $ \Box $    
   

\bigskip
\bigskip

\bibliographystyle{unsrt.bst}
\bibliography{biblio}

\end{document}